\documentclass[11pt]{amsart}
\usepackage{amsmath}
\usepackage{amsthm}
\usepackage{amssymb}
\usepackage{amsfonts}
\usepackage{amsxtra}  
\usepackage{epsfig}
\usepackage{verbatim}

\theoremstyle{plain}
\newtheorem{theorem}[equation]{Theorem}
\newtheorem{proposition}[equation]{Proposition}
\newtheorem{lemma}[equation]{Lemma}

\newtheorem{definition}[equation]{Definition}
\theoremstyle{remark}
\newtheorem*{remark}{Remark}
\numberwithin{equation}{section}

\include{header}

\newcommand{\Tnorm}{|\!|\!|}
\newcommand{\dbar}{\bar{\partial}}
\newcommand{\dbarstar}{\bar{\partial}^{\star}}
\newcommand{\dbarstarp}{\bar{\partial}_{\phi}^{\star}}
\newcommand{\dbarstarpd}{\bar{\partial}_{\phi_{\delta}}^{\star}}
\newcommand{\dbarstarpe}{\bar{\partial}_{\phi_{\eta}}^{\star}}
\newcommand{\phidel}{\phi_{\delta}}

\newcommand{\lamdel}{\lambda_{\delta}}
\newcommand{\phieta}{\phi_{\eta}}
\newcommand{\esdel}{S_{\delta}}
\newcommand{\Omegaa}{\Omega_{a}}
\newcommand{\zetadel}{\zeta_{\delta}}
\newcommand{\kernel}{\text{ker}\,}
\newcommand{\ssubset}{\subset\subset}
\newcommand{\psum}{\sideset{}{^{\prime}}{\sum}}

\begin{document}

\bibliographystyle{plain}       

\title[$\dbar$-Neumann Operator]{A Sufficient Condition for Subellipticity of the $\dbar$-Neumann Operator}
\author{Anne-Katrin Herbig}
\address{Department of Mathematics, University of Michigan, Ann Arbor, Michigan 48109}
\email{herbig@umich.edu}
\date{}
\begin{abstract}
We give a sufficient condition for subelliptic estimates for the $\dbar$-Neumann operator on smoothly bounded, pseudoconvex domains in $\mathbb{C}^n$. This condition is a quantified version of McNeal's condition $(\tilde{P})$ for compactness of the $\dbar$-Neumann operator, and it extends Catlin's sufficiency condition for subellipticity as it is less stringent.
\end{abstract}
\maketitle

\section{Introduction} 
 Let $\Omega\ssubset\mathbb{C}^{n}$ be a smoothly bounded domain. Suppose that $p\in b\Omega$ is a point in the boundary of $\Omega$, and that $b\Omega$ is pseudoconvex near $p$. We shall show that the existence of a certain family of functions near the boundary point $p$ implies that a subelliptic estimate for the $\dbar$-Neumann operator holds near that point.

 The $\dbar$-Neumann operator $N_{p,q}$ is the inverse of the complex Laplacian $\dbar\dbarstar+
 \dbarstar\dbar$ for $(p,q)$-forms. 
  Establishing
  the existence of the $\dbar$-Neumann operator 
  leads to a particular solution of the Cauchy-Riemann equations, but 
  just in the $L^{2}$-sense. Thus one is not just interested in the 
  existence of such an 
  $L^{2}$-solution $u$ for given data $f$, but one is 
  also interested in the kind of regularity statements that can be made about 
  $u$ when $f$ 
  is regular; for notation and details on the $\dbar$-Neumann problem see section \ref{S:sectionPrelim}.

  On domains with certain geometric conditions on the 
  boundary, the question of existence of a solution to the 
  $\dbar$-Neumann problem was settled through the works of 
  H\"ormander \cite{Hor65}, Kohn \cite{Koh63,Koh64} and Morrey 
  \cite{Mor58}. In fact, H\"ormander's results in \cite{Hor65} imply 
  that there exists a bounded operator $N_{p,q}$ on 
  $L_{p,q}^{2}(\Omega)$, which inverts the complex Laplacian 
  under the assumption that $\Omega$ is a bounded, 
  pseudoconvex domain.
  
  In the following, we will be concerned only with the local 
  regularity question for the $\dbar$-Neumann problem, i.e. 
  conditions on $\Omega$ which imply that $u:=N_{p,q}f$ is 
  smooth wherever $f$ is. A fundamental step concerning this question 
  was done by Kohn and Nirenberg. They showed in \cite{Koh-Nir} that, 
  if a so-called \emph{subelliptic estimate of order $\epsilon$} holds for 
  the $\dbar$-Neumann problem on a 
  neighborhood $V$ of a given point $p$ in $b\Omega$, then 
  $f_{|_{V}}\in H^{s}_{p,q}(V)$ implies $N_{p,q}f_{|_{V'}}\in 
  H^{s+2\epsilon}_{p,q}(V')$ for $V'\ssubset V$; here $H^{s}_{p,q}$ 
  denotes the $L^{2}$-Sobolev space of order $s$ on $(p,q)$-forms.
  Thus it is natural to inquire about subelliptic estimates for the 
  $\dbar$-Neumann problem.
  
  Denote by $\mathcal{D}^{p,q}(V\cap\bar{\Omega})$ 
  the set of smooth $(p,q)$-forms 
  $u$, which are supported in $V\cap\bar{\Omega}$, such that
  $u$ belongs to the domain of $\dbarstar$. A subelliptic estimate of order
  $\epsilon>0$  near $p\in b\Omega$ is said to hold, if
  \begin{eqnarray}\label{D:SubEllE}
    \Tnorm u\Tnorm_{\epsilon}^{2}\leq C(\|\dbar u\|^{2}+\|\dbarstar u\|^{2})    
    \,\,\text{ for all }u\in\mathcal{D}^{p,q}(V\cap\bar{\Omega}),
  \end{eqnarray}
  where the norm on the left hand side
   is the tangential $L^{2}$-Sobolev norm of order $\epsilon$.
  
   The most general result concerning subelliptic estimates for the 
  $\dbar$-Neumann problem was obtained by Catlin \cite{Cat87}.
  He showed that the existence of a certain, uniformly bounded family 
  of functions  $\{\lamdel\}$ on 
  a pseudoconvex domain is sufficient for a subelliptic estimate to 
  hold.
 Moreover, Catlin proved that one can construct such a family of 
 functions on any smoothly bounded, pseudoconvex domain, which is of 
  finite type in the sense of D'Angelo \cite{D'An82}.

  We extend Catlin's sufficiency result by replacing the boundedness condition on the weight functions 
  $\lamdel$ with that of \emph{self-bounded complex gradient}, a weaker 
  condition which allows unbounded families of functions. This notion 
  was introduced by McNeal in \cite{McN02}.
  \begin{definition}
    Let $\Omega\ssubset\mathbb{C}^{n}$ be a smoothly bounded domain. A 
    plurisubharmonic function $\phi\in C^{2}(\Omega)$ is said to have 
    a self-bounded complex gradient, if there exists a constant $C>0$ 
    such that
    \begin{eqnarray}\label{D:SbcG}
      |\sum_{k=1}^{n}\frac{\partial\phi}{\partial z_{k}}(z)\xi_{k}|^{2}
      \leq C\sum_{k,l=1}^{n}\frac{\partial^{2}\phi}{\partial z_{k}
      \partial\bar{z}_{l}}(z)\xi_{k}\bar{\xi}_{l}
    \end{eqnarray}	
    holds for all $\xi\in\mathbb{C}^{n}$, $z\in\Omega$.
  We write $|\partial\phi|_{i\partial\dbar\phi}\leq \sqrt{C}$ 
  when we mean (\ref{D:SbcG}).
  \end{definition}
  Notice that, if $\lambda\in C^{2}(\Omega)$ is plurisubharmonic and 
  bounded, then $\phi=e^{\lambda}$ has a self-bounded complex gradient 
  with $C=\sup_{z\in\Omega} e^{\lambda(z)}$. Furthermore, notice the 
  behavior of inequality (\ref{D:SbcG}) 
  under scaling; replacing $\phi$ by $t\phi$ for $t>0$, the left 
  hand side of (\ref{D:SbcG}) is quadratic in $t$, while the right 
  hand side is linear in $t$.
  
  The main result in this paper is the following:
  \begin{theorem}\label{T:MainTheorem}
   Let $\Omega\ssubset\mathbb{C}^{n}$ be a smoothly bounded domain. 
    Let $p$ be a given point in $b\Omega$ and suppose that 
    $b\Omega\cap U$ is pseudoconvex, where $U$ is a neighborhood of 
    $p$. Denote by $\esdel$ the set $\lbrace z\in\Omega\,|\,-\delta<
    r(z)<0\rbrace$, where $r$ is a fixed, smooth defining function of $\Omega$. 
    Assume that for all $\delta>0$ sufficiently 
    small there exists a plurisubharmonic function
    $\phidel\in C^{2}(\bar{\Omega}\cap U)$, such 
    that
    \begin{enumerate}
      \item[(i)] $|\partial\phidel|_{i\partial\dbar\phidel}^{2}\leq C$, where the constant 
      $C>0$ is independent of $\delta$,
      \item[(ii)] for all smooth $(p,q)$-forms u, $z\in \esdel\cap 
        U$ and for some $\epsilon\in(0,\frac{1}{2}{]}$
	\begin{eqnarray*}
	  \psum_{|I|=p,|J|=q-1}\sum_{k,l=1}^{n}
	  \frac{\partial^{2}\phidel}{\partial z_{k}\partial\bar{z}_{l}}(z)
	  u_{I,kJ}\bar{u}_{I,lJ}
	  \geq c\delta^{-2\epsilon}|u|^{2},
	\end{eqnarray*}  
	  where the  constant $c>0$ does not depend on $\delta$ or $u$.    
    \end{enumerate}	
    Then there exists a neighborhood $V\ssubset U$ of $p$ such that
    a subelliptic estimate of order $\epsilon$ holds.
  \end{theorem}    
  The only difference between Theorem \ref{T:MainTheorem} and Catlin's sufficiency result is that we  substituted the uniform boundedness condition on $\{\lamdel\}$ by condition (i). 
 The existence of Catlin's family of functions 
  $\lbrace\lamdel\rbrace$ implies the existence of the above family 
  $\lbrace\phidel\rbrace$ by setting $\phidel=e^{\lamdel}$. One reason, however, to generalize the Theorem of Catlin is to establish sharper subelliptic estimates in various geometric situations.
  
The uniform boundedness of $\{\lamdel\}$ is crucial for Catlin's proof as it lets him transform estimates  
with weights of the form $e^{-\lamdel}$ into unweighted estimates.  Families of functions which have a self-bounded complex gradient are in general not uniformly bounded, and so Catlin's proof does not work. However,  McNeal found a duality argument in \cite{McN02}, which allows one to pass to unweighted estimates  from estimates with weights, when the weight functions have a self-bounded complex gradient.
    
The paper is structured as follows. In section 
\ref{S:sectionPrelim} we review briefly the setting of the $\dbar$-Neumann 
problem. 
In section \ref{S:sectionWeightinq} we derive two 
weighted $L^{2}$-inequalities, which are specific for weights having 
a self-bounded complex gradient. Using those inequalities we obtain  
two versions of compactness estimates on $\dbarstar N_{q}$ and 
$\dbarstar N_{q+1}$ in section \ref{S:sectiondbarstarN}.
In section \ref{S:sectionDest} we convert these compactness 
estimates to a family of $L^{2}$-estimates in terms of the Dirichlet 
form. 
With those estimates at hand 
we complete the proof of Theorem \ref{T:MainTheorem} in 
section \ref{S:sectionSubest}. In the last section we consider an example domain to see how the functions $\{\phidel\}$ can be constructed.

I am deeply indebted to  J.D. McNeal for his support and encouragement. I have enjoyed and greatly  benefitted from our discussions during the last years.

\section{Preliminaries}\label{S:sectionPrelim}
  Let $\Omega\ssubset\mathbb{C}^{n}$ be a smoothly bounded domain, 
  i.e. $\Omega$ is bounded and there is a smooth function $r$ such that
  $\Omega=\lbrace z\in\mathbb{C}^{n}\,|\,r(z)<0\rbrace$ and $\nabla r\neq 
  0$ whenever $r=0$. \\
  Let $0\leq p,q\leq n$. We write an arbitrary $(p,q)$-form $u$ as
  \begin{eqnarray}\label{N:pqformnotation}
    u=\psum_{|I|=p,|J|=q}u_{I,J}dz^{I}\wedge d\bar{z}^{J},  
  \end{eqnarray}
  where $I=\lbrace i_{1},\dots,i_{p}\rbrace$, $J=\lbrace 
  j_{1}\dots,j_{q}\rbrace$ and 
  $dz^{I}=dz^{i_{1}}\wedge\dots\wedge dz^{i_{p}}$, 
  $d\bar{z}^{J}=d\bar{z}^{j_{1}}\wedge\dots\wedge d\bar{z}^{j_{q}}$. 
  Here $\sum^{\prime}$ means that we only sum over strictly increasing 
  index sets. We define the coefficients $u_{I,J}$ for arbitrary index sets 
  $I$ and $J$, so that the $u_{I,J}$'s are antisymmetric functions of $I$ 
  and $J$.\\
  Let $\Lambda^{p,q}(\bar{\Omega})$ and $\Lambda_{c}^{p,q}(\Omega)$
denote the $(p,q)$-forms with coefficients in
  $C^{\infty}(\bar{\Omega})$ and $C_{c}^{\infty}(\Omega)$, 
  respectively. We use the 
  pointwise inner product $\langle\,.,.\,\rangle$ defined by
   $ \langle dz^{k},dz^{l}\rangle=\delta^{k}_{l}
    =\langle d\bar{z}^{k},d\bar{z}^{l}\rangle$.  
  By linearity we extend this inner product to $(p,q)$-forms.
  The global $L^{2}$-inner product on $\Omega$ is defined by
  \begin{eqnarray*}
      (u,v)_{\Omega}=\int_{\Omega}\langle u,v\rangle dV,
  \end{eqnarray*}    
  where $dV$ is the euclidean volume form. 
  The $L^{2}$-norm of a $u\in\Lambda_{c}^{p,q}(\Omega)$ on $\Omega$ is then given by 
  $\|u\|_{\Omega}^{2}=(u,u)_{\Omega}$ and we define $L_{p,q}^2(\Omega)$ to be the completion of $\Lambda_{c}^{p,q}(\Omega)$ under the $L^2$-norm; we drop the subscript $\Omega$, 
  when there is no reason for confusion.\\
  If $\phi\in C^{2}(\bar{\Omega})$, we denote by 
  $L_{p,q}^{2}(\Omega,\phi)$ the space of $(p,q)$-forms $u$ such that
  \begin{eqnarray*}
    \|u\|_{\phi,\Omega}^{2}=(u,u)_{\phi,\Omega}:=
    \|u e^{-\frac{\phi}{2}}\|_{\Omega}^{2}=
    \int_{\Omega}\langle u,u\rangle e^{-\phi}dV<\infty.
  \end{eqnarray*}
  Notice that the weighted $L^{2}$-space, $L_{p,q}^{2}(\Omega,\phi)$, equals 
  $L_{p,q}^{2}(\Omega)$.\\
  Let $u\in\Lambda^{p,q}(\bar{\Omega})$, then the $\dbar$-operator is 
  defined as
  \begin{eqnarray*}
      \dbar_{p,q}u=\dbar u:=\psum_{|I|=p,|J|=q}\sum_{k=1}^{n}
      \bar{\partial}_{k}u_{I,J}\,d\bar{z}^{k}\wedge dz^{I}\wedge 
      d\bar{z}^{J},
  \end{eqnarray*}  
  where $\bar{\partial}_{k}:=\frac{\partial}{\partial\bar{z}_{k}}$, and 
  $u$ is expressed as in (\ref{N:pqformnotation}). 
  Observe that $\dbar^{2}=0$.
  We extend the differential operator $\dbar$, still denoted by 
  $\dbar$, to act on non-smooth forms in the sense of distributions. 
  Then, by restricting the domain of $\dbar$ to those forms $g\in 
  L_{p,q}^{2}(\Omega)$, where $\dbar g$ in the distributional sense
  belongs to $L_{p,q+1}^{2}(\Omega)$, $\dbar$ becomes an operator on Hilbert 
  spaces at each form level. Note that $\dbar$ is a densely defined 
  operator on $L_{p,q}^{2}(\Omega)$, since the compactly supported 
  forms $\Lambda_{c}^{p,q}(\Omega)$ are in $\text{Dom}(\dbar)$. 
  Moreover, $\dbar$ is a closed operator, because differentiation is 
  a continuous map in the distributional sense.\\
  Thus we can define the Hilbert space adjoint, $\dbarstar$, to 
  $\dbar$ with respect to the $L^{2}$-inner product on the 
  appropriate form level in the usual way:
  
  we say that $u\in L_{p,q+1}^{2}(\Omega)$ belongs to the domain of 
  $\dbarstar$, i.e. $u\in\text{Dom}(\dbarstar)$, 
  if there exists a constant $C>0$ so that
  \begin{eqnarray}\label{E:defDomdbarstar}
    |(\dbar w,u)|\leq 
    C\|w\| \text{ holds for all } w\in\text{Dom}(\dbar).
  \end{eqnarray}  
  By the Riesz representation theorem it follows, that, if 
  $u\in\text{Dom}(\dbarstar)$, there exists a unique $v\in 
  L_{p,q}^{2}(\Omega)$, such that
  \begin{eqnarray*}
      (w,v)=(\dbar w, u)
  \end{eqnarray*}    
  holds for all $w\in\text{Dom}(\dbar)$; we write $\dbarstar u$ for 
  $v$. 
  This
  reveals that 
  certain boundary conditions must hold on any smooth $(p,q+1)$-form, which 
  belongs to $\text{Dom}(\dbarstar)$. In fact, one can show that
    $u\in\mathcal{D}^{p,q+1}(\Omega):=\text{Dom}(\dbarstar)\cap
    \Lambda^{p,q+1}(\bar{\Omega})$ holds if and only if
    \begin{eqnarray*}
      \sum_{k=1}^{n}u_{I,kJ}\frac{\partial r}{\partial 
      z_{k}}=0\,\text{on }b\Omega	
    \end{eqnarray*}	
    for all $I$ and $J$ which are strictly increasing index sets of 
    length $p$ and $q$, respectively. Here, $r$ is a defining function of $\Omega$.
 
 The Hilbert space adjoint, $\dbarstarp$, 
  to $\dbar$ with respect to the $L^{2}(\Omega,\phi)$-inner product is defined by
  $\dbarstarp=e^{\phi}\dbarstar e^{-\phi}$.
  In view of (\ref{E:defDomdbarstar}) it is easy to see that
  $\text{Dom}(\dbarstar)=\text{Dom}(\dbarstarp)$ holds.\\
    
  Now we are ready to formulate the $\dbar$-Neumann problem. It is the 
  following: given $f\in L_{p,q}^{2}(\Omega)$, find $u\in 
  L_{p,q}^{2}(\Omega)$ such that the following holds
  \begin{eqnarray}\label{E:dbarNeumannproblem}
    \begin{cases}
      (\dbar\dbarstar+\dbarstar\dbar)u=f\\
      u\in\text{Dom}(\dbar)\cap\text{Dom}(\dbarstar)\\
      \dbar u\in\text{Dom}(\dbarstar),\,\dbarstar u\in\text{Dom}(\dbar)
    \end{cases}	
  \end{eqnarray} 
  The complex Laplacian, $\square_{p,q}:=\dbar\dbarstar+\dbarstar\dbar$, 
  is itself elliptic, but the boundary conditions, which are implied 
  by membership to $\text{Dom}(\dbarstar)$, are not. 
  The ellipticity of $\square_{p,q}$ implies that G$\mathring{\text{a}}$rding's inequality holds in the 
  interior of $\Omega$, i.e.
   \begin{eqnarray}\label{E:Garding}
     \|u\|_{1}^{2}\lesssim\|\dbar u\|^{2}+\|\dbarstar u\|^{2}
    \text{ for } u\in\Lambda_{c}^{p,q}(\Omega),
   \end{eqnarray}
  where $\|.\|_{1}$ denotes the usual $L^{2}$-Sobolev $1$-norm. 
  We remark, though, (\ref{E:Garding}) does not 
  hold for general $u\in\mathcal{D}^{p,q}(\Omega)$.  However, a substitute estimate, (\ref{E:SubEll}) below, does hold for $u\in\mathcal{D}^{p,q}(\Omega)$.
  
  Let $p\in b\Omega$. We may choose 
 a neighborhood $U$ of $p$ and a local coordinate system 
  $(x_{1},\ldots,x_{2n-1},r)\in\mathbb{R}^{2n-1}\times\mathbb{R}$, 
  such that the last coordinate is a local defining function of the 
  boundary. Call $(U,(x,r))$ a special boundary chart. We shall 
  denote the dual variable of $x$ by $\xi$, and define $\langle 
  x,\xi\rangle:=\sum_{j=1}^{2n-1}x_{j}\xi_{j}$.
  For $f\in C_{c}^{\infty}(U\cap\bar{\Omega})$ we define the 
  tangential Fourier transform of $f$ by
  \begin{eqnarray*}
    \tilde{f}(\xi,r):=\int_{\mathbb{R}^{2n-1}}e^{-2\pi i\langle 
    x,\xi\rangle}f(x,r)dx.  
  \end{eqnarray*}  
  Via the tangential Bessel potential $\Lambda_{t}^{s}$ of order $s$,
  \begin{eqnarray*}
      (\Lambda_{t}^{s}f)(x,r):=
      \int_{\mathbb{R}^{2n-1}}e^{2\pi i\langle x,\xi\rangle}
      (1+|\xi|^{2})^{\frac{s}{2}}\tilde{f}(\xi,r)d\xi,
  \end{eqnarray*}
  we can define the tangential $L^{2}$-Sobolev norm of $f$ of order 
  $s$ by
  \begin{eqnarray*}
    \Tnorm f\Tnorm_{s}^{2}:=\|\Lambda_{t}^{s}f\|^{2}
    =\int_{-\infty}^{0}\int_{\mathbb{R}^{2n-1}}
    (1+|\xi|^{2})^{s}|\tilde{f}(\xi,r)|^{2}d\xi dr.
  \end{eqnarray*} 
  A subelliptic estimate of order $\epsilon>0$ holds if there exists $C>0$ such that
  \begin{eqnarray}\label{E:SubEll}
  \Tnorm u\Tnorm_{\epsilon}^{2}\leq C \|\dbar u\|^2+\|\dbarstar u\|^2
  \end{eqnarray}
  for $u\in\mathcal{D}^{p,q}(\Omega)$ supported near the boundary point $p$.

  From here on, we restrict our considerations to $(0,q)$-forms. The 
  system (\ref{E:dbarNeumannproblem}) does not see the $dz$'s and the 
  general case for $(p,q)$-forms can be derived easily.
  For notational ease we shall
  write $u_{J}$, instead of $u_{0,J}$, for the components of a 
  $(0,q)$-from $u$. We shall denote the Dirichlet form 
  associated to $\square_{0,q}$ as usual by $Q(.,.)$, i.e. $Q(u,v):=
  (\dbar u,\dbar v)+(\dbarstar u,\dbarstar v)$ for 
  $u,v\in\mathcal{D}^{0,q}(\Omega)$.\\ 
  For quantities $A$ and $B$ we use 
  the notation $|A|\lesssim|B|$ to mean $|A|\leq C|B|$ for some 
  constant $C>0$, which is independent of relevant parameters. It will 
  be specifically mentioned or clear from the context, what those 
  parameters are.
  Furthermore, we call the elementary inequality
   $ |AB|\leq\eta A^{2}+\frac{1}{4\eta}B^{2}$ for 
   $ \eta>0$ the (sc)-(lc) inequality.

\section{Basic estimates}\label{S:sectionWeightinq}
  In this section, we derive two basic weighted inequalities for forms 
  in $\mathcal{D}^{0,q}(\Omega)$. 
  We will make extensive use of these inequalities in 
  our proof of subellipticity. 
 Our starting point is the following Proposition 
 \ref{P:BasicEst}, which has been derived by McNeal in \cite{McN02}.
  \begin{proposition}\label{P:BasicEst}
    Let $\Omega\ssubset\mathbb{C}^{n}$ be a smoothly bounded, pseudoconvex 
    domain, and suppose that $\phi\in C^{2}(\bar{\Omega})\cap 
    PSH(\Omega)$. 
    If $|\partial\phi|_{i\partial\dbar\phi}\leq 1$, then	
      \begin{eqnarray}
          \frac{1}{2}\psum_{|I|=q-1}\int_{\Omega} 
          \sum_{k,l=1}^{n}\frac{\partial^{2}\phi}{\partial 
          z_{k}\partial\bar{z}_{l}} u_{kI} \bar{u}_{lI}e^{-2\phi}dV
          \leq \|\dbar u\|_{2\phi}^{2}
          +3\|\dbar_{\phi}^{\star}u\|_{2\phi}^{2} \label{F:SbgE1}
        \end{eqnarray}
	holds for all $u\in\mathcal{D}^{0,q}(\Omega)$.
      \end{proposition}

  We remark that inequality (\ref{F:SbgE1}) is one of the key points 
  leading to the subelliptic estimate. In fact, this inequality will 
  be used in section \ref{S:sectiondbarstarN} enabling us to obtain 
  ``good'' estimates near the boundary. In the following,  
  we derive a
  G$\mathring{\text{a}}$rding-like weighted inequality. This 
  inequality is also crucial as it will give us ``good'' 
  estimates in the interior.
  \begin{proposition}\label{P:IntEst}
    Let $\Omega\ssubset\mathbb{C}^{n}$ be a smoothly bounded, pseudoconvex 
    domain, and suppose that $\phi\in C^{2}(\bar{\Omega})\cap 
    PSH(\Omega)$ satisfies $|\partial\phi|_{i\partial\dbar\phi}\leq 
    \frac{1}{\sqrt{24}}$. Then for all 
    $u\in\Lambda_{c}^{0,q}(\Omega)$, 
    it holds that
    \begin{eqnarray}\label{E:IntEst}
      \|u e^{-\phi}\|_{1}^{2}\lesssim \|\dbar u\|_{2\phi}^{2}
      +\|\dbarstarp u\|_{2\phi}^{2},
    \end{eqnarray}
    where $\|.\|_{1}$ denotes the $L^{2}$-Sobolev $1$-norm on $\Omega$.
  \end{proposition}
  For the proof of Proposition \ref{P:IntEst} we need to introduce the 
  Hodge-Star Operator $\star$, that is the map
    \begin{displaymath}
	\star:\Lambda^{p,q}(\bar{\Omega})\longrightarrow\Lambda^{n-p,n-q}
	(\bar{\Omega})
    \end{displaymath}
  defined by $\psi\wedge\star\varphi=\langle\psi,\varphi\rangle dV$ 
  for $\psi,\varphi\in\Lambda^{p,q}(\bar{\Omega})$. The basic 
  properties of the Hodge-Star Operator are summarized in the 
  following lemma.
  \begin{lemma}\label{L:Hodge}
    \begin{enumerate}
      \item[(i)] $\star\star=(-1)^{p+q}\text{ id}$ on 
        $\Lambda^{p,q}(\bar{\Omega})$, 
      \item[(ii)] $|\varphi|=|\star\varphi|$ for 
        $\varphi\in\Lambda^{p,q}(\bar{\Omega})$, 
        where $|\varphi|^2=\langle\varphi,\varphi\rangle$,
      \item[(iii)] $ \dbarstar=-\star\dbar\star$ 
      on $\Lambda_{c}^{p,q}(\bar{\Omega})$.	
    \end{enumerate}	
  \end{lemma}    
  A proof of Lemma \ref{L:Hodge} can  be found in  \cite{Che-Sha01}, chapter 9.
  \begin{proof}[Proof of Proposition \ref{P:IntEst}]
    Let $u\in\Lambda_{c}^{0,q}(\bar{\Omega})$. 
    By G$\mathring{\text{a}}$rding's inequality (\ref{E:Garding}),
    we have
    \begin{eqnarray*}
      \|u e^{-\phi}\|_{1}^{2}\lesssim 
      \|\dbar (u e^{-\phi})\|^{2}
      +\|\dbarstar (u e^{-\phi})\|^{2}
      =\|\dbar (u e^{-\phi})\|^{2}
      +\|\dbarstarp u\|_{2\phi}^{2}.
    \end{eqnarray*}
    Thus we just need to consider the term $\|\dbar (u e^{-\phi})\|^{2}$.
    For that define $v\in\Lambda_{c}^{n,n-q}(\bar{\Omega})$ by $v=\star 
    u$. Here we denote the coefficients of $v$ by $v_{J}$ for $|J|=n-q$. 
    Then, by Lemma \ref{L:Hodge} and commuting, it follows
    \begin{eqnarray*}
      \|\dbar (u e^{-\phi})\|^{2}& =&
      \|\dbarstar (v e^{-\phi})\|^{2}
      \lesssim \|\dbarstar v\|_{2\phi}^{2}+
      \|\lbrack\dbarstar ,\phi\rbrack v\|_{2\phi}^{2}\\
      &=&\|-\star\dbar \star v\|_{2\phi}^{2}
       +\psum_{|J|=n-q-1}
      \|\sum_{l=1}^{n}\frac{\partial\phi}{\partial z_{l}} 
      v_{lJ}\|_{2\phi}^{2}\\
      &\leq&\|\dbar u\|_{2\phi}^{2}+\psum_{|J|=n-q-1}
      \int_{\Omega}\sum_{k,l=1}^{n}\frac{\partial^{2}\phi}
      {\partial z_{k}\partial\bar{z}_{l}} v_{kJ}\bar{v}_{lJ} e^{-2\phi} 
      dV, 
    \end{eqnarray*}
    where the last step follows from $\phi$ having a self-bounded 
    complex gradient. Note that $v\in\mathcal{D}^{n,n-q}(\Omega)$, 
    since $v$ is identically zero on the boundary of $\Omega$. Hence 
    we can apply inequality (\ref{F:SbgE1}):
    \begin{eqnarray*}
      \psum_{|J|=n-q-1}
      \int_{\Omega}\sum_{k,l=1}^{n}\frac{\partial^{2}\phi}
      {\partial z_{k}\partial\bar{z}_{l}} v_{kJ}\bar{v}_{lJ} e^{-2\phi} 
      dV
      \leq
      2\|\dbar v\|_{2\phi}^{2}+6\|\dbarstarp  
      v\|_{2\phi}^{2}.
    \end{eqnarray*}	
    Since $|\partial\phi|_{i\partial\dbar\phi}\leq\frac{1}{\sqrt{24}}$, 
    it follows that
    \begin{eqnarray*}
      \|\dbarstarp v\|_{2\phi}^{2}
      &\leq&
      2\|\dbarstar v\|_{2\phi}^{2}+2\|\lbrack\dbarstar,\phi\rbrack 
      v\|_{2\phi}^{2}\\
      &=&
      2\|\dbarstar v\|_{2\phi}^{2}+2\psum_{|J|=n-q-1}\|\sum_{l=1}^{n}
      \frac{\partial\phi}{\partial z_{l}}v_{lJ}\|_{2\phi}^{2}\\
      &\leq&
      2\|\dbarstar v\|_{2\phi}^{2}+
      \frac{1}{12}\psum_{|J|=n-q-1}
      \int_{\Omega}\sum_{k,l=1}^{n}\frac{\partial^{2}\phi}
      {\partial z_{k}\partial\bar{z}_{l}} v_{kJ}\bar{v}_{lJ} e^{-2\phi} 
      dV.
    \end{eqnarray*}
    Thus we obtain
    \begin{eqnarray*}
      \psum_{|J|=n-q-1}
      \int_{\Omega}\sum_{k,l=1}^{n}\frac{\partial^{2}\phi}
      {\partial z_{k}\partial\bar{z}_{l}} v_{kJ}\bar{v}_{lJ} e^{-2\phi} 
      dV
      &\leq&
      4\|\dbar v\|_{2\phi}^{2}+24\|\dbarstar  
      v\|_{2\phi}^{2}\\
      &=&
      4\|\dbarstar u\|_{2\phi}^{2}+24\|\dbar u\|_{2\phi}^{2}¥
    \end{eqnarray*}	
    where the second line holds by Lemma \ref{L:Hodge}. So we are 
    left with estimating the term 
    $\|\dbarstar u\|_{2\phi}^{2}$. As before, we just need to commute:
    \begin{eqnarray*}
      \|\dbarstar u\|_{2\phi}^{2}
      &\lesssim&
      \|\dbarstarp u\|_{2\phi}^{2}
      +\|\lbrack\dbarstar ,\phi\rbrack u\|_{2\phi}^{2}
      =
      \|\dbarstarp u\|_{2\phi}^{2}
      +\psum_{|I|=q-1}\|\sum_{l=1}^{n}\frac{\partial\phi}
      {\partial z_{l}}u_{lI}\|_{2\phi}^{2}\\
      &\leq&
      \|\dbarstarp u\|_{2\phi}^{2}
      +\psum_{|I|=q-1}\int_{\Omega}\sum_{k,l=1}^{n}
      \frac{\partial^{2}\phi}{\partial z_{k}\partial\bar{z}_{l}}
      u_{kI}\bar{u}_{lI} e^{-2\phi} dV,
    \end{eqnarray*}	
    which, again, follows by the self-bounded complex gradient 
    condition of $\phi$. To finish we use inequality (\ref{F:SbgE1}) 
    again, that is
    \begin{eqnarray*}
      \psum_{|I|=q-1}\int_{\Omega}\sum_{k,l=1}^{n}
      \frac{\partial^{2}\phi}{\partial z_{k}\partial\bar{z}_{l}}
      u_{kI}\bar{u}_{lI} e^{-2\phi} dV
      \lesssim
      \|\dbar u\|_{2\phi}^{2}
      +\|\dbarstarp u\|_{2\phi}^{2}.
    \end{eqnarray*}	
    Collecting all our estimates, we obtain
    \begin{displaymath}
      \|u e^{-\phi}\|_{1}^{2}\lesssim \|\dbar u\|_{2\phi}^{2}
      +\|\dbarstarp u\|_{2\phi}^{2}\;
    \text{ for }\; u\in\Lambda_{c}^{0,q}(\bar{\Omega}).
    \end{displaymath}
  \end{proof}    
  Since the $L^{2}$-Sobolev $1$-norm dominates the $L^{2}$-norm, 
  (\ref{E:IntEst}) implies that
  \begin{eqnarray*}
      \|ue^{-\phi}\|^{2}\lesssim\|\dbar 
      u\|_{2\phi}^{2}+\|\dbarstarp u\|_{2\phi}^{2}
  \end{eqnarray*}    
  holds for all $u\in\Lambda_{c}^{0,q}(\Omega)$. In the following, we 
  show that this inequality is in fact true for all 
  $u\in\mathcal{D}^{0,q}(\Omega)$. 
  \begin{proposition}
    Let $\Omega\ssubset\mathbb{C}^{n}$ be a smoothly bounded, pseudoconvex 
    domain, and suppose that $\phi\in C^{2}(\bar{\Omega})\cap 
    PSH(\Omega)$ satisfies $|\partial\phi|_{i\partial\dbar\phi}\leq
    \frac{1}{\sqrt{2}}$. Then 
    for $u\in\mathcal{D}^{0,q}(\Omega)$ it holds that
    \begin{eqnarray}
      \|u\|_{2\phi}^{2} \lesssim \|\dbar u\|_{2\phi}^{2}
      +\|\dbarstarp u\|_{2\phi}^{2} \label{F:SbgE3}.
    \end{eqnarray}	  
  \end{proposition}    
  \begin{proof}
     Set $\psi_{t}(z)=\phi(z)+t|z|^{2}$ for $t>0$. Then $\psi_{t}$ is 
     strictly plurisubharmonic, since for
     $\xi\in\mathbb{C}^{n}$, $z\in\Omega$ it holds
     \begin{eqnarray*}
       \sum_{k,l=1}^{n}\frac{\partial^{2}\psi_{t}}{\partial z_{k}
       \partial\bar{z}_{l}}(z)\xi_{k}\bar{\xi}_{l}
       =\sum_{k,l=1}^{n}\frac{\partial^{2}\phi}{\partial z_{k}
       \partial\bar{z}_{l}}(z)\xi_{k}\bar{\xi}_{l}
       +t|\xi|^{2}.
     \end{eqnarray*} 
     Moreover, we observe that
     \begin{eqnarray*}
       |\sum_{k=1}^{n}\frac{\partial\psi_{t}}
       {\partial z_{k}}(z)\xi_{k}|^{2}
       \leq 2|\sum_{k=1}^{n}\frac{\partial\phi}{\partial 
       z_{k}}(z)\xi_{k}|^{2}
       +2t^{2}|z|^{2}|\xi|^{2}.
     \end{eqnarray*}	 
     Since $\Omega$ 
     is a bounded domain, we can choose a $t>0$, such that 
      $24t|z|^2\leq 1$ holds for all $z\in\Omega$. 
     Then $|\partial\psi_{t}|_{i\partial\dbar\psi_{t}}\leq 1$, and thus
     inequality (\ref{F:SbgE1}) holds for $\psi_{t}$. That is
     \begin{eqnarray*}
       \frac{1}{2}\psum_{|I|=q-1}\int_{\Omega} 
          \sum_{k,l=1}^{n}\frac{\partial^{2}\psi_{t}}{\partial 
          z_{k}\partial\bar{z}_{l}} u_{kI} \bar{u}_{lI}e^{-2\psi_{t}}dV
          \leq \|\dbar u\|_{2\psi_{t}}^{2}
          +3\|\dbar_{\psi_{t}}^{\star}u\|_{2\psi_{t}}^{2}	 
     \end{eqnarray*}
     Note that $e^{-2t|z|^{2}}$ is bounded from above by $1$ and that $\phi$ is 
     plurisubharmonic on $\Omega$. Hence it follows that
     \begin{eqnarray*}
       \frac{1}{2}\int_{\Omega} 
          t|u|^{2}e^{-2\psi_{t}}dV
          &\leq& \|\dbar u\|_{2\phi}^{2}
          +3\|\dbar_{\psi_{t}}^{\star}u\|_{2\psi_{t}}^{2}\\
	  &\leq&
	  \|\dbar u\|_{2\phi}^{2}+6\|\dbarstarp u\|_{2\phi}^{2}
	  +6\|\lbrack\dbarstar,(t|z|^{2})\rbrack u\|_{2\psi_{t}}^{2}.
     \end{eqnarray*}
     By our choice of $t$ we can estimate the last term
     \begin{eqnarray*}
       6\|\lbrack\dbarstar,(t|z|^{2})\rbrack u\|_{2\psi_{t}}^{2}
       =6\psum_{|I|=q-1}\|\sum_{k=1}^{n}\frac{\partial(t|z|^{2})}
       {\partial z_{k}} u_{kI}\|_{2\psi_{t}}^{2}
       \leq
       \frac{1}{4}t\|u\|_{2\psi_{t}}^{2}.
     \end{eqnarray*}	
     Therefore it holds that
     \begin{eqnarray*}
       \frac{1}{4}\int_{\Omega}t|u|^{2}e^{-2\psi_{t}}
       \leq
       \|\dbar u\|_{2\phi}^{2}+6\|\dbarstarp u\|_{2\phi}^{2}
     \end{eqnarray*}	 
     Since $e^{-t|z|^{2}}$ is bounded from 
     below on $\Omega$, our claim follows.  
 \end{proof}    
  
\section{Estimates for $\dbarstar N_q$}\label{S:sectiondbarstarN}
  By a compactness estimate for $\dbarstar N_{q}$ we mean the 
  following: for all $\eta>0$ there exists a $C(\eta)>0$ such that
  \begin{eqnarray}\label{D:CompEst}
    \|\dbarstar N_{q}\alpha\|\lesssim\eta\|\alpha\|+C(\eta)
    \|\alpha\|_{-s}
  \end{eqnarray} 
  for all $\alpha\in L_{0,q}^{2}(\Omega)$. Here $\|.\|_{-s}$, $s>0$, denotes 
  the $L^{2}$-Sobolev norm of order $-s$. The constant in 
  $\lesssim$ does depend on $s$ but not on $\alpha$, $\eta$ or $C(\eta)$.
  The family of estimates (\ref{D:CompEst}) is equivalent to 
  $\dbarstar N_{q}$ being a compact operator from 
  $L_{0,q}^{2}(\Omega)$ to $L_{0,q-1}^{2}(\Omega)$; for a proof see 
  for instance \cite{McN02}. We remark that for compactness of 
  $\dbarstar N_{q}$ it is sufficient to establish (\ref{D:CompEst}) 
  for $\dbar$-closed forms $\alpha\in L_{0,q}^{2}(\Omega)$, see \cite{McN02}.

  In this section, we derive with the aid of our weighted estimates 
  from section 
  \ref{S:sectionWeightinq} two versions of compactness estimates 
  for $\dbarstar N_{q}$. We 
  start out with a quantified version of (\ref{D:CompEst}), i.e. we 
  describe $C(\eta)$ for each $\eta$.\\
  
  Since the weight functions $\lbrace\phidel\rbrace$ are just defined 
  on $\Omega\cap U$, where $U$ is a neighborhood of a given  
  $p\in b\Omega$ (see hypotheses in Theorem \ref{T:MainTheorem}), we need 
  to restrict our considerations to an approximating 
  subdomain of $\Omega$, which lies in $U$.
  \begin{proposition}\label{P:ApprSubD}
    Suppose that $\Omega\subset\mathbb{C}^{n}$ is a smoothly bounded domain. 
    Let $p$ be a point in $b\Omega$ and suppose that $b\Omega\cap U$ 
    is pseudoconvex, where $U$ is a neighborhood of $p$.
    Then there exists a smoothly bounded, pseudoconvex domain
    $\Omegaa\subset\Omega\cap U$ with $\Omegaa\ssubset U$ satisfying the following 
    properties
    \begin{enumerate}
      \item[(1)]
      $b\Omega\cap b\Omegaa$ contains a neighborhood of $p$
      in $b\Omega$,
      \item[(2)]
      all points in $b\Omegaa\setminus b\Omega$ are strongly 
      pseudoconvex.
    \end{enumerate}
  \end{proposition}
  A proof of Proposition \ref{P:ApprSubD} can be found in \cite{McN92}.
  We call such a domain $\Omegaa$ an approximating subdomain 
  associated to $(\Omega,p,U)$. The crucial feature, for our current purposes, 
  of such an approximating subdomain $\Omegaa$ is that 
  it is a smoothly bounded, pseudoconvex domain. 
  Therefore we can apply
  the inequalities 
  (\ref{F:SbgE1}), (\ref{E:IntEst}) and (\ref{F:SbgE3})
  on $\Omegaa$ using the 
  $\phidel$'s as weight functions. We remark that for using these 
  inequalities a rescaling of the 
  $\phidel$'s might be necessary, so that 
  $|\partial\phidel|_{i\partial\dbar\phidel}\leq\frac{1}{\sqrt{24}}$ 
  holds for all $\delta>0$ sufficiently small. 
  \begin{theorem}\label{T:dbarstarNE}
    Assume the hypotheses of Theorem \ref{T:MainTheorem}.
    Let $\Omegaa$ be an approximating subdomain associated to 
    $(\Omega, p,U)$. Then there exists a neighborhood 
    $V\ssubset U$ of $p$, such that for 
    $\alpha\in L_{0,q}^{2}(\Omegaa)$, $\dbar$-closed 
    and supported in $V\cap\bar{\Omega}_{a}$, the following estimate 
    holds:
    \begin{eqnarray}
      \|\dbarstar N^{\Omega_{a}}_{q}\alpha\|_{\Omega_{a}}^{2}
      \lesssim
      \delta^{2\epsilon}\|\alpha\|_{\Omega_{a}}^{2}
      +\delta^{-2+2\epsilon}\|\alpha\|_{-1,\Omega_{a}}^{2}.
    \end{eqnarray}
    The constant in $\lesssim$ neither depends on $\alpha$ nor 
    $\delta$.
  \end{theorem}
  \begin{proof}
    For notational ease we shall write $\|.\|$ for 
    $\|.\|_{\Omega_{a}}$ and $N_{q}$ for $N_{q}^{\Omega_{a}}$. 
    Let $W\ssubset U$ be a neighborhood of $p$, such that 
    $W\cap\Omega\subset\Omegaa$ and 
    $\overline{W}\cap b\Omegaa\ssubset b\Omega$. Also, let $V\ssubset W$ 
    be a neighborhood of $p$ and $\alpha\in L_{0,q}^{2}(\Omegaa)$ 
    be a $\dbar$-closed form, 
    which is supported in $V\cap\bar{\Omega}_{a}$.   
    Define the functional $F : (\lbrace e^{-\frac{\phidel}{2}}
    \dbarstarpd u\,|\,u\in\mathcal{D}^{0,q}(\Omegaa)\rbrace, \|.\|_{\phidel})
    \longrightarrow\mathbb{C}$ by 
    \begin{displaymath}
      F(e^{-\frac{\phidel}{2}}\dbarstarpd u)=(u,\alpha)_{\phidel}.
    \end{displaymath}  
    We start with showing that $F$ satisfies the 
    following estimate
    \begin{eqnarray}\label{E:Fbound}
      |F(e^{-\frac{\phidel}{2}}\dbarstarpd u)|
      \lesssim
      \|e^{-\frac{\phidel}{2}}\dbarstarpd u\|_{\phidel}
      (\delta^{\epsilon}\|\alpha\|
      +\delta^{-1+\epsilon}\|\alpha\|_{-1}).
    \end{eqnarray}
    Recall that $S_{\delta}=\lbrace 
    z\in\Omegaa\,|\,-\delta<r(z)<0\rbrace$, where $r$ is the fixed defining 
    function of $\Omega$.
    Let $\chi\in C_{c}^{\infty}(W)$ such that $\chi\equiv 1$ on $V$ 
    and $\chi\geq 0$. 
    Recall that the support of $\alpha$ is in $V$. Then, by the 
    generalized Cauchy-Schwarz inequality, we obtain
    \begin{eqnarray*}
      |F(e^{-\frac{\phidel}{2}}\dbarstarpd u)|
      &=&|(u,\alpha)_{\phidel}|\\
      &\lesssim&\|ue^{-\phidel}\|^{W\cap\esdel}\|\alpha\|
      +\|e^{-\phidel}\chi 
      u\|_{1}^{\Omegaa\backslash\esdel}\|\alpha\|_{-1}.
    \end{eqnarray*}
    In view of our claim (\ref{E:Fbound}) we need to estimate the 
    terms $\|ue^{-\phidel}\|^{W\cap\esdel}$ and $\|e^{-\phidel}\chi 
      u\|_{1}^{\Omegaa\backslash\esdel}$ appropriately.\\
      1. Estimating $\|ue^{-\phidel}\|^{W\cap\esdel}$:  
      Recall that $\phidel$ has a self-bounded complex gradient on 
      $\Omegaa\subset U\cap\Omega$ by 
      hypothesis (i). Hence inequality (\ref{F:SbgE1}) holds, and 
       the plurisubharmonicity of $\phidel$ implies then, that
      \begin{displaymath}
	\psum_{|I|=q-1}\int_{W\cap\esdel}\sum_{k,l=1}^{n}
	\frac{\partial^{2}\phidel}{\partial z_{k}\partial\bar{z}_{l}}
	u_{kI}\bar{u}_{lI}e^{-2\phidel}dV
        \lesssim \|\dbar u\|_{2\phidel}^{2}
        +\|\dbarstarpd u\|_{2\phidel}^{2}  
      \end{displaymath}
      holds uniformly for all $\delta>0$ small.
      Invoking hypothesis (ii) and noting that $W\subset U$ yields
      \begin{eqnarray}\label{F:ColardelE}
	 \|u\|_{2\phidel}^{W\cap\esdel}
	 \lesssim \delta^{\epsilon}(\|\dbar u\|_{2\phidel}
        +\|\dbarstarpd u\|_{2\phidel}).
      \end{eqnarray}
      2. Estimating $\|e^{-\phidel}\chi 
      u\|_{1}^{\Omegaa\backslash\esdel}$: Let 
      $h_{\delta}:\mathbb{R}_{0}^{+}\longrightarrow\lbrack 0,1\rbrack$
      be a smooth function with $h_{\delta}(x)=0$ for $x\in\lbrack 
      0,\frac{\delta}{2}\rbrack$ and 
      $h_{\delta}(x)=1$ for $x\geq\delta$. We can choose $h_{\delta}$ 
      such that $|h_{\delta}'|\lesssim\delta^{-1}$. Define
      $\zetadel\in C^{\infty}(\bar{\Omega}_{a})$ by 
      $\zetadel(z)=h_{\delta}(-r(z))$, where $r$ is the fixed defining 
      function of $\Omega$. Note that
      \begin{eqnarray}\label{F:ZetaDerE}
	\left|\frac{\partial\zetadel}{\partial x_{j}}\right|
	\lesssim
	\delta^{-1}\left|\frac{\partial r}{\partial x_{j}}\right|
	\lesssim\delta^{-1}
      \end{eqnarray}
      holds on $\Omegaa$ for all $j\in\lbrace 1,\dots,2n\rbrace$.
    Clearly, we have
    \begin{eqnarray*}
      \|e^{-\phidel}\chi 
      u\|_{1}^{\Omegaa\backslash\esdel}
      \leq\|e^{-\phidel}\zetadel\chi 
      u\|_{1}.
    \end{eqnarray*}
    Since $\zetadel\cdot \chi$ is identically zero near the boundary of $\Omegaa$, 
    we can use our weighted G$\mathring{\text{a}}$rding's inequality 
    (\ref{E:IntEst}) to start estimating
    \begin{eqnarray*}
      \|e^{-\phidel}\zetadel\chi u\|_{1}^{2}
      &\lesssim&
      \|\dbar(\zetadel\chi u)\|_{2\phidel}^{2}
      +\|\dbarstarpd(\zetadel\chi u)\|_{2\phidel}^{2}\\
      &\lesssim&
      \|\dbar u\|_{2\phidel}^{2}+\|\dbarstarpd u\|_{2\phidel}^{2}
      +\sum_{j=1}^{2n}\left(\|\frac{\partial\zetadel}{\partial 
      x_{j}}\chi u\|_{2\phidel}^{2}+\|\frac{\partial\chi}{\partial 
      x_{j}}u\|_{2\phidel}^{2}\right)\\
      &\lesssim&\|\dbar u\|_{2\phidel}^{2}+\|\dbarstarpd u\|_{2\phidel}^{2}
      +\sum_{j=1}^{2n}\max_{z\in\bar{\Omega}_{a}}\left|
      \frac{\partial\zetadel}{\partial x_{j}}\right|^{2}
      (\|u\|_{2\phidel}^{W\cap\esdel})^{2}
      +\|u\|_{2\phi}^{2}.
    \end{eqnarray*}
    The last estimate holds since $\chi$ is supported in $W$ and 
    $\frac{\partial\zetadel}{\partial x_{j}}=0$ on $\Omegaa\backslash
    \esdel$. By  the inequalities (\ref{F:SbgE3}) and (\ref{F:ZetaDerE}),
    it follows 
    \begin{eqnarray*}
      \|e^{-\phidel}\zetadel\chi 
      u\|_{1}^{2}
      \lesssim
      \|\dbar u\|_{2\phidel}^{2}+\|\dbarstarpd u\|_{2\phidel}^{2}
      +\delta^{-2}(\|u\|_{2\phidel}^{W\cap\esdel})^{2}
    \end{eqnarray*}
    for all $\delta>0$ small enough.
    Using the estimate (\ref{F:ColardelE}) for 
    $\|u\|_{2\phidel}^{W\cap\esdel}$, we obtain
    \begin{eqnarray*}
      \|e^{-\phidel}\zetadel\chi 
      u\|_{1}^{2}
      \lesssim
      \delta^{-2+2\epsilon}(\|\dbar u\|_{2\phidel}^{2}+
      \|\dbarstarpd u\|_{2\phidel}^{2}),	
    \end{eqnarray*}
    thus we can conclude 
    \begin{eqnarray}\label{F:OhneColardelE}
      (\|e^{-\phidel}\chi 
      u\|_{1}^{\Omegaa\backslash\esdel})^{2}
      \lesssim
      \delta^{-2+2\epsilon}(\|\dbar u\|_{2\phidel}^{2}+
      \|\dbarstarpd u\|_{2\phidel}^{2}).	
    \end{eqnarray}	
    Write $u=u_{1}+u_{2}$, where $u_{1}\in \kernel\dbar$ and 
    $u_{2}\perp_{\phidel} \kernel\dbar$. Note that 
    $u_{1}\in\mathcal{D}^{0,q}(\Omegaa)$. Thus, since $\alpha\in
    \kernel\dbar$, we get, using the
    estimates (\ref{F:ColardelE}) and (\ref{F:OhneColardelE}),
    \begin{eqnarray*}
      |(u,\alpha)_{\phidel}|
      =
      |(u_{1},\alpha)_{\phidel}|
      \lesssim
      \|\dbarstarpd u_{1}\|_{2\phidel}
      (\delta^{\epsilon}\|\alpha\|+\delta^{-1+\epsilon}\|\alpha\|_{-1}).
    \end{eqnarray*}
    However, $u_{2}\perp_{\phidel} \kernel\dbar$, therefore we get 
    $\|\dbarstarpd u\|_{2\phidel}^{2}=\|\dbarstarpd u_{1}\|_{2\phidel}^{2}$.
     Hence our claimed inequality (\ref{E:Fbound}) holds:
    \begin{eqnarray*}
      |F(e^{-\frac{\phidel}{2}}\dbarstarpd u)|=|(u,\alpha)_{\phidel}|
      \lesssim
      \|e^{-\frac{\phidel}{2}}\dbarstarpd u\|_{\phidel}
      (\delta^{\epsilon}\|\alpha\|
      +\delta^{-1+\epsilon}\|\alpha\|_{-1}).
    \end{eqnarray*}	
    That is, $F$ is a bounded linear functional on 
    $(\lbrace e^{-\frac{\phidel}{2}}
    \dbarstarpd u\,|\,u\in\mathcal{D}^{0,q}(\Omegaa)\rbrace, 
    \|.\|_{\phidel})$, which is a subset of 
    $L_{0,q-1}^{2}(\Omegaa,\phidel)$. By the Hahn-Banach Theorem, $F$ extends 
    to a bounded linear functional on 
    $L_{0,q-1}^{2}(\Omegaa,\phidel)$ with the same bound. The Riesz 
    representation theorem yields, that there exists a unique $v\in
    L_{0,q-1}^{2}(\Omegaa,\phidel)$ such that for all $g\in
    L_{0,q-1}^{2}(\Omegaa,\phidel)$
    \begin{eqnarray*}
      F(g)&=&(g,v)_{\phidel},\\
      \|v\|_{\phidel}^{2}&\lesssim&\delta^{2\epsilon}\|\alpha\|^{2}
      +\delta^{-2+2\epsilon}\|\alpha\|_{-1}^{2}.
    \end{eqnarray*}
    In particular, we get for all $u\in\mathcal{D}^{0,q}(\Omegaa)$
    \begin{eqnarray*}
      (u,\dbar(e^{-\frac{\phidel}{2}}v))_{\phidel}
      =(e^{-\frac{\phidel}{2}}\dbarstarpd u,v)_{\phidel}
      =(u,\alpha)_{\phidel}.	
    \end{eqnarray*}
    Note that $\mathcal{D}^{0,q}(\Omegaa)$ is dense in $L_{(0,q)}^{2}
    (\Omegaa,\phidel)$. Hence,
    setting $s=e^{-\frac{\phidel}{2}}v$, it follows that $\dbar 
    s=\alpha$ in the distributional sense and
    \begin{eqnarray*}
      \|s\|^{2}\lesssim
      \delta^{2\epsilon}\|\alpha\|^{2}
      +\delta^{-2+2\epsilon}\|\alpha\|_{-1}^{2}.
    \end{eqnarray*}
    But the minimal $L^{2}(\Omegaa)$-solution, $\dbarstar N_{q}\alpha$, to 
    the $\dbar$-problem for $\alpha$ on $\Omegaa$ 
    must also satisfy this estimate; 
    that is
    \begin{eqnarray}\label{E:locdbarstarNqE}
      \|\dbarstar N_{q}\alpha\|^{2}\lesssim
      \delta^{2\epsilon}\|\alpha\|^{2}
      +\delta^{-2+2\epsilon}\|\alpha\|_{-1}^{2}.
    \end{eqnarray}
  \end{proof}
  \begin{remark}
    Observe that the only point where the form level $q$ of
    the $(0,q)$-forms comes into play, is in 
     hypothesis (ii) of Theorem \ref{T:MainTheorem}.  
    Notice that 
    this condition on the complex hessian of $\phidel$ near the 
    boundary also holds for $(0,q+1)$-forms. 
    Thus
    by a  proof analogous to the above, we obtain the following: there exists a 
    neighborhood $V\ssubset U$ of $p$ such that for all $\beta\in
    L_{0,q+1}(\Omegaa)$, which are $\dbar$-closed and supported in 
    $V\cap\bar{\Omega}_{a}$, the following estimate holds
    \begin{eqnarray}\label{E:locdbarstarNq+1E}
      \|\dbarstar N_{q+1}^{\Omega_{a}}\beta\|^{2}
      \lesssim
      \delta^{2\epsilon}\|\beta\|_{\Omega_{a}}^{2}
      +\delta^{-2+2\epsilon}\|\beta\|_{-1,\Omega_{a}}^{2}.
    \end{eqnarray}	
  \end{remark}  
  
  These families of estimates, (\ref{E:locdbarstarNqE}) and 
  (\ref{E:locdbarstarNq+1E}), are the heart of the matter for our 
  proof of subellipticity. 
  But to convert these estimates on $\dbarstar N_{q}^{\Omega_{a}}$ and 
  $\dbarstar N_{q+1}^{\Omega_{a}}$ to usable estimates on 
  $\mathcal{D}^{0,q}(\Omega)$, we shall need  
  exact regularity of the operator $\dbarstar\dbar N_{q}^{\Omega_{a}}$.
  By exact regularity we mean that $\dbarstar\dbar N_{q}^{\Omegaa}$ 
  preserves the 
  $L^{2}$-Sobolev spaces. 
    
  Kohn showed in \cite{Koh84}, that exact 
  regularity of $\dbarstar\dbar N_{q}^{\Omega}$ follows from compactness of 
  $N_{q}^{\Omega}$ on $L_{0,q}^{2}(\Omega)$, if $\Omega$ is a 
  smoothly bounded, pseudoconvex domain. It is an easy consequence of 
  the formula
  \begin{eqnarray*}
    N_{q}=(\dbar N_{q-1})(\dbarstar N_{q})+(\dbarstar 
    N_{q+1})(\dbar N_{q}),   
  \end{eqnarray*}    
  that compactness of the operators  $\dbarstar N_{q}$ and $\dbarstar 
  N_{q+1}$ implies compactness of $N_{q}$.

  The estimates (\ref{E:locdbarstarNqE}) 
  and (\ref{E:locdbarstarNq+1E}) do not imply compactness as they do not hold for all $\dbar$-closed           
  forms in $L_{0,q}^{2}(\Omegaa)$ and $L_{0,q+1}^{2}(\Omegaa)$, respectively. 
 However, we show below that 
  $N_{q}^{\Omega_{a}}$ is a compact operator on 
  $L_{0,q}^{2}(\Omegaa)$ by using a proof similar to the one of 
  Theorem \ref{T:dbarstarNE}. The crucial property of the 
  approximating subdomain $\Omega_{a}$ for this argument is that 
  $\Omegaa$ is strongly pseudoconvex off the boundary of $\Omega$. In 
  particular, we use Kohn's result that near a point in the boundary of 
  strong pseudoconvexity a subelliptic estimate of order 
  $\frac{1}{2}$ holds.
  \begin{proposition}\label{P:LocNComp}
    Assume that the hypotheses of Theorem \ref{T:MainTheorem} hold. Let
    $\Omegaa$ be an approximating subdomain associated to 
    $(\Omega,p,U)$.
    Then the $\dbar$-Neumann operator $N_{q}^{\Omegaa}$ 
    is a compact operator on $L_{0,q}^{2}(\Omegaa)$.
  \end{proposition} 
  \begin{proof}
    As before, we write $N_{q}$ for $N_{q}^{\Omegaa}$, and $\|.\|$ 
    for $\|.\|_{\Omega_{a}}$.
    We start out with showing that $\dbarstar N_{q}$ is a compact 
    operator. 
    By the remark following (\ref{D:CompEst}) we obtain compactness of 
    $\dbarstar N_{q}$, if we can show that for all 
    $\eta>0$ there exists a $C(\eta)>0$ such that 
    \begin{eqnarray*}
      \|\dbarstar N_{q}\alpha\|\lesssim\eta\|\alpha\|+
      C(\eta)\|\alpha\|_{-\frac{1}{2}}
    \end{eqnarray*}	
    holds for all $\dbar$-closed $\alpha\in L^{2}_{0,q}(\Omegaa)$.
    
    Let $\eta>0$ be given. By our hypotheses there exists a function $\phieta
    \in C^{2}(\bar{\Omega}_{a})\cap PSH(\Omegaa)$ which has a 
    self-bounded complex gradient and satisfies
    \begin{eqnarray}\label{E:etablowup}
      \psum_{|I|=q-1}\sum_{k,l=1}^{n}
      \frac{\partial^{2}\phieta}{\partial z_{k}
      \partial\bar{z}_{l}}(z)u_{kI}\bar{u}_{lI}
      \geq\eta^{-2}|u|^{2} \text{  for }u\in\Lambda^{0,q}(\Omegaa)
    \end{eqnarray}	  
    on a strip $S_{\eta'}=\lbrace z\in\Omegaa\cap\Omega\,|\,
      -\eta'<r(z)<0\rbrace$ for some $\eta'>0$ chosen small enough, 
      depending on $\eta$. Here $r$ is the fixed defining function of $\Omega$.
      
    Let $\alpha$ be a $\dbar$-closed $(0,q)$-form with coefficients in 
    $L^{2}(\Omegaa)$. Define the linear functional
    $F:(\lbrace e^{-\frac{\phieta}{2}}\dbarstarpe u\,
    |\,u\in\mathcal{D}^{0,q}(\Omegaa)\rbrace,\|\,.\,\|_{\phieta})
    \longrightarrow\mathbb{C}$ by
    \begin{eqnarray*}
      F(e^{-\frac{\phieta}{2}}\dbarstarpe u)=(u,\alpha)_{\phieta}.	
    \end{eqnarray*}
    We shall show that $F$ is a bounded functional satisfying
    \begin{eqnarray}\label{E:schlFbound}
      |F(e^{-\frac{\phieta}{2}}\dbarstarpe u)|
      \lesssim
      \|e^{-\frac{\phieta}{2}}\dbarstarpe u\|_{\phieta}
      (\eta\|\alpha\|+C(\eta)\|\alpha\|_{-\frac{1}{2}})
    \end{eqnarray}	
    for some $C(\eta)>0$.
    For that let $\chi\in C^{\infty}(\Omegaa)$ be a non-negative 
    function such that 
    $\chi=1$ on $\Omegaa\setminus S_{\eta'}$ and $\chi=0$ on
    $S_{\frac{\eta'}{2}}$. Then
    \begin{eqnarray*}
      |F(e^{-\frac{\phieta}{2}}\dbarstarpe u)|
      &=&
      |(u,\alpha)^{S_{\eta'}}|+|(u,\alpha)^{\Omegaa\setminus 
      S_{\eta'}}|\\
      &\lesssim&
      \|ue^{-\phieta}\|^{S_{\eta'}}\|\alpha\|+
      \|\chi u e^{-\phieta}\|_{\frac{1}{2}}
      \|\alpha\|_{-\frac{1}{2}},
    \end{eqnarray*}	
    where the second line follows by the generalized Cauchy-Schwarz 
    inequality. In view of our claimed inequality 
    (\ref{E:schlFbound}), we need to get control of the terms 
    $\|ue^{-\phieta}\|^{S_{\eta'}}$ and 
    $\|\chi u e^{-\phieta}\|_{\frac{1}{2}}$. 
    
    Since $\phieta\in C^{2}(\bar{\Omega}_{a})\cap PSH(\Omegaa)$ 
    has a self-bounded complex gradient and $\Omegaa$ is 
    pseudoconvex, we can use inequality (\ref{F:SbgE1}) to 
    estimate $\|ue^{-\phieta}\|^{S_{\eta'}}$: 
    \begin{eqnarray*}
      \psum_{|I|=q-1}\int_{S_{\eta'}}
      \sum_{k,l=1}^{n}\frac{\partial^{2}\phieta}{\partial z_{k}
      \partial \bar{z}_{l}}u_{kI}\bar{u}_{lI} e^{-2\phieta} dV
      \lesssim
      \|\dbar u\|_{2\phieta}^{2}+\|\dbarstarpe u\|_{2\phieta}^{2}. 	
    \end{eqnarray*}	
    By inequality (\ref{E:etablowup}) it follows 
    \begin{eqnarray*}
	(\|u e^{-\phieta}\|^{S_{\eta'}})^{2}
	\lesssim
	\eta(\|\dbar u\|_{2\phieta}^{2}+\|\dbarstarpe u\|_{2\phieta}^{2}).
    \end{eqnarray*}
    
    In order to estimate $\|\chi u e^{-\phieta}\|_{\frac{1}{2}}$, note that 
    $\text{supp}\chi\cap b\Omegaa\ssubset b\Omegaa\setminus
    b\Omega$ and recall that, by our choice of $\Omegaa$, we 
    have that $b\Omegaa\backslash b\Omega$ is strongly pseudoconvex. 
    Thus an subelliptic estimate of order $\frac{1}{2}$ holds for 
    $\chi u e^{-\phieta}$:
    \begin{eqnarray*}
      \|\chi u e^{-\phieta}\|_{\frac{1}{2}}^{2}
      &\lesssim&
      \|\dbar(\chi u e^{-\phieta})\|^{2}+\|\dbarstar(\chi u 
      e^{-\phieta})\|^{2}\\
      &\lesssim&
      C(\eta)^{2}(\|\dbar u\|_{2\phieta}^{2}
      +\|\dbarstarpe u\|_{2\phieta}^{2}
      +\|u\|_{2\phieta}^{2})\\
      &\lesssim&
      C(\eta)^2
      (\|\dbar u\|_{2\phieta}^{2}+\|\dbarstarpe u\|_{2\phieta}^{2}),
    \end{eqnarray*}	
    where the last line follows by inequality (\ref{F:SbgE3}).
    
    Now we are set up for proving inequality (\ref{E:schlFbound}). 
    Write $u=u_{1}+u_{2}$, where $u_{1}\in\kernel\dbar$ and 
    $u_{2}\perp_{\phieta}\kernel\dbar$. Thus, since $\alpha\in\kernel
    \dbar$, we get, using our above estimates for the terms 
    $\|u e^{-\phieta}\|^{S_{\eta'}}$ and 
    $\|\chi u e^{-\phieta}\|_{\frac{1}{2}}$,
    \begin{eqnarray*}
      |F(e^{-\frac{\phieta}{2}}\dbarstarpe u)|
      =|(u_{1},\alpha)_{\phieta}|
      &\lesssim&\|u_{1} e^{-\phieta}\|^{S_{\eta'}}\|\alpha\|
      +\|\chi u_{1} e^{-\phieta}\|_{\frac{1}{2}}\|\alpha\|_{-\frac{1}{2}}\\
      &\lesssim&
      \|\dbarstarpe u_{1}\|_{2\phieta}(\eta\|\alpha\|
      +C(\eta)\|\alpha\|_{-\frac{1}{2}}).
    \end{eqnarray*}
    Recall that $\|\dbarstarpe u\|_{2\phieta}
    =\|\dbarstarpe u_{1}\|_{2\phieta}$ holds, since 
    $u_{2}\perp_{\phieta}\kernel\dbar$.
    This implies our claimed inequality (\ref{E:schlFbound}).   
    By arguments analogous to the ones in the  proof of Theorem \ref{T:dbarstarNE} it follows that
      \begin{eqnarray*}
      \|\dbarstar N_{q}\alpha\|\lesssim
      \eta\|\alpha\|+C(\eta)\|\alpha\|_{-\frac{1}{2}}.
    \end{eqnarray*}
    holds for all $\dbar$-closed forms $\alpha\in L_{0,q}^{2}(\Omegaa)$.
    Thus $\dbarstar N_{q}$ is a compact 
    operator from $L_{0,q}^{2}(\Omegaa)$ to $L_{0,q-1}^{2}(\Omegaa)$. 
    A similar proof yields the 
    compactness of $\dbarstar N_{q+1}$. Therefore $N_{q}$, 
    the $\dbar$-Neumann operator on 
    $\Omegaa$, is a compact operator on $L_{0,q}^{2}(\Omegaa)$. 
  \end{proof}

\section{Estimates on $\mathcal{D}^{0,q}(\Omega)$}\label{S:sectionDest}
  In this section we convert the families of estimates, 
  (\ref{E:locdbarstarNqE}) and (\ref{E:locdbarstarNq+1E}), obtained 
  in section \ref{S:sectiondbarstarN} to estimates 
  for forms in $\mathcal{D}^{0,q}(\Omega)$. As already mentioned in 
  section \ref{S:sectiondbarstarN}, we need exact regularity 
  to hold for operators related to $N_{q}^{\Omegaa}$.  We begin with 
  a result of Kohn.
  \begin{proposition}\label{L:Kohnsglobreg}
    Suppose $\Omega\ssubset\mathbb{C}^{n}$ is a smoothly bounded, 
    pseudoconvex domain, such that its $\dbar$-Neumann operator, 
    $N_{q}$, is compact on $L_{0,q}^{2}(\Omega)$. Let $s>0$, then the 
    following holds
    \begin{enumerate} 
      \item[(1)] if $\beta\in H_{0,q}^{s}(\Omega)$, then
                 $\|\dbarstar\dbar N_{q}\beta\|_{s}\lesssim\|\beta\|_{s}$,
      \item[(2)] if $\beta\in H_{0,q-1}^{s}(\Omega)$, then
                 $\|N_{q}\dbar\beta\|_{s}\lesssim\|\beta\|_{s}$.
    \end{enumerate}
    Here, the constants in $\lesssim$ depend on $s$ but not on 
    $\beta$. 
  \end{proposition}
  A proof of Proposition \ref{L:Kohnsglobreg} is contained in 
  \cite{Koh84}. An easy consequence of Proposition \ref{L:Kohnsglobreg} 
  is the exact regularity of the $L^{2}$-adjoint operators of $\dbarstar\dbar 
  N_{q}$ and $N_{q}\dbar$ in the $L^{2}$-Sobolev spaces of 
  negative order. In particular, the following holds.
  \begin{lemma}\label{L:negexRN}
    Suppose $\Omega\ssubset\mathbb{C}^{n}$ is a smoothly bounded, 
    pseudoconvex domain, such that its $\dbar$-Neumann operator, 
    $N_{q}$, is compact on $L_{0,q}^{2}(\Omega)$. Then, if $\alpha
    \in \Lambda^{0,q}(\bar{\Omega})$, it follows that
    \begin{eqnarray}
      \|\dbarstar N_{q}\alpha\|_{-1} 
      &\lesssim&\|\alpha\|_{-1},\label{E:dbarstarN-1}\\
      \|\dbar\dbarstar N_{q}\alpha\|_{-1}
      &\lesssim&\|\alpha\|_{-1}.\label{E:dbardbarstarN-1}
    \end{eqnarray}	
  \end{lemma}  
  \begin{proof}
      Let $\alpha\in\Lambda^{0,q}(\bar{\Omega})$. Then
      \begin{eqnarray*}
	\|\dbarstar N_{q}\alpha\|_{-1}
	=\sup\lbrace (\dbarstar N_{q}\alpha,\beta)\,|\,
	\beta\in H_{0,q-1}^{1}(\Omega),\|\beta\|_{1}\leq 1\rbrace.
      \end{eqnarray*}
      Since $\beta\in H_{0,q-1}^{1}(\Omega)$ is in $\text{Dom}(\dbar)$, we 
      obtain
      \begin{eqnarray*}
	(\dbarstar N_{q}\alpha,\beta)
	=(N_{q}\alpha,\dbar\beta)
	=(\alpha,N_{q}\dbar\beta)
	\lesssim\|\alpha\|_{-1}\|N_{q}\dbar\beta\|_{1},
      \end{eqnarray*}
      by the generalized Cauchy-Schwarz 
      inequality. Proposition \ref{L:Kohnsglobreg}, part (2),
      yields exact regularity for $N_{q}\dbar$, in particular 
      $\|N_{q}\dbar\beta\|_{1}\lesssim\|\beta\|_{1}$ holds
        for all $\beta\in H_{0,q-1}^{1}(\Omega)$. Thus we have
      \begin{eqnarray*}
	\|\dbarstar N_{q}\alpha\|_{-1}
	\lesssim
	\sup\lbrace \|\alpha\|_{-1}\|\beta\|_{1}\,|\,
	\beta\in H_{0,q-1}^{1}(\Omega),\|\beta\|_{1}\leq 1\rbrace
	=\|\alpha\|_{-1},
      \end{eqnarray*}	  
     which proves (\ref{E:dbarstarN-1}).\\
     The proof of (\ref{E:dbardbarstarN-1}) is very similar. 
     Since $\alpha=(\dbar\dbarstar N_{q}+\dbarstar\dbar 
     N_{q})\alpha$, it holds that
     \begin{eqnarray*}
       \|\dbar\dbarstar N_{q}\alpha\|_{-1}
       =\|\alpha-\dbarstar\dbar N_{q}\alpha\|_{-1}
       \leq\|\alpha\|_{-1}+\|\dbarstar\dbar N_{q}\alpha\|_{-1},
     \end{eqnarray*}
     where
     \begin{eqnarray*}
       \|\dbarstar\dbar N_{q}\alpha\|_{-1}
       =\sup\lbrace(\dbarstar\dbar N_{q}\alpha,\beta)\,|\,\beta
       \in H_{0,q}^{1}(\Omega),\|\beta\|_{1}\leq 1\rbrace.
     \end{eqnarray*}	 
     As before, note that $\beta\in H_{0,q}^{1}(\Omega)$ is in
     $\text{Dom}(\dbar)$. Moreover, since 
     $\alpha\in\Lambda^{0,q}(\bar{\Omega})$, it holds that 
     $\dbar N_{q}\alpha=N_{q+1}\dbar\alpha$. Thus we obtain
     \begin{eqnarray*}
       (\dbarstar\dbar N_{q}\alpha,\beta)
       &=&
       (N_{q+1}\dbar\alpha,\dbar\beta)
       =(\alpha,\dbarstar N_{q+1}\dbar\beta)
       =(\alpha,\dbarstar\dbar N_{q}\beta)\\
       &\lesssim&
       \|\alpha\|_{-1}\|\dbarstar\dbar N_{q}\beta\|_{1}.
     \end{eqnarray*}
     Part (1) of Proposition \ref{L:Kohnsglobreg} 
     tells us that
      $ \|\dbarstar\dbar N_{q}\beta\|_{1}\lesssim\|\beta\|_{1}$	 
     holds for all $\beta\in H_{0,q}^{1}(\Omega)$.
     Hence it follows 
     \begin{eqnarray*}
       \|\dbarstar\dbar N_{q}\alpha\|_{-1}
       \lesssim\sup\lbrace\|\alpha\|_{-1}\|\beta\|_{1}\,|\,\beta
       \in H_{0,q}^{1}(\Omega),\|\beta\|_{1}\leq 1\rbrace
       =\|\alpha\|_{-1},
     \end{eqnarray*}
     which proves (\ref{E:dbardbarstarN-1}).
  \end{proof}
  Recall that we showed in Proposition \ref{P:LocNComp} that the 
  $\dbar$-Neumann operator, $N_{q}^{\Omegaa}$, associated to the 
  approximating subdomain $\Omegaa$ is compact. Therefore, the exact 
  regularity results (\ref{E:dbarstarN-1}) and 
  (\ref{E:dbardbarstarN-1}) hold for $N_{q}^{\Omegaa}$. Now we are 
  ready to derive estimates for forms in $\mathcal{D}^{0,q}(\Omega)$.
  \begin{proposition}\label{P:NpassingQ}
    Assume the hypotheses of Theorem \ref{T:MainTheorem}.
    Then there exists a neighborhood $W\ssubset U$ of 
    $p$, such that for all sufficiently small $\delta>0$ and $\eta>0$
    \begin{eqnarray*}
      \|u\|_{\Omega}^{2}
      \lesssim
      \frac{\delta^{2\epsilon}}{\eta}(\|\dbar u\|_{\Omega}^{2}
      +\|\dbarstar u\|_{\Omega}^{2}
      +\delta^{-2}\|\dbar u\|_{-1,\Omega}^{2})
      +\eta\delta^{-2}\|u\|_{-1,\Omega}^{2}
    \end{eqnarray*}
    holds for $u\in\mathcal{D}^{0,q}(\Omega)$ supported in 
    $W\cap\bar{\Omega}$. Here, the constant in $\lesssim$ does depend on $\eta$ but not on $\delta$.
  \end{proposition}
  \begin{proof}
    Recall that
    Theorem \ref{T:dbarstarNE} and the following remark say that if $\Omegaa$ is an 
    approximating subdomain associated to $(\Omega,p,U)$, then 
    there exists a neighborhood 
    $V\ssubset U$ of $p$ such that
    \begin{eqnarray}
      \|\dbarstar N_{q}^{\Omegaa}\alpha\|_{\Omegaa}^{2}&\lesssim&
      \delta^{2\epsilon}\|\alpha\|_{\Omegaa}^{2}
      +\delta^{-2+2\epsilon}\|\alpha\|_{-1,\Omegaa}^{2},\label{E:dbarstarNq}\\
      \|\dbarstar N_{q+1}^{\Omegaa}\beta\|_{\Omegaa}^{2}&\lesssim&
      \delta^{2\epsilon}\|\beta\|_{\Omegaa}^{2}
      +\delta^{-2+2\epsilon}\|\beta\|_{-1,\Omegaa}^{2} \label{E:dbarstarNq+1}
    \end{eqnarray}
    hold for all $\alpha\in L_{0,q}^{2}(\Omegaa)$ and $\beta\in 
    L_{0,q+1}^{2}(\Omegaa)$, which are $\dbar$-closed and supported in
    $V\cap\bar{\Omega}_{a}$.
    For notational ease we denote the $L^{2}$-norm on $\Omegaa$ by 
    $\|.\|$ and write $N_{q}$ for the $\dbar$-Neumann operator on 
    $\Omegaa$.
    
    Recall that $V$ in Theorem \ref{T:dbarstarNE}
    was chosen such that $V\cap 
    b\Omegaa\ssubset b\Omega$. 
    Let $W\ssubset V$ be a neighborhood of $p$, and 
    $\zeta\in C_{c}^{\infty}(V)$, $\zeta\geq 0$ and $\zeta\equiv 1$ on $W$. Let 
    $u\in\mathcal{D}^{0,q}(\Omega)$ be supported in 
    $W\cap\bar{\Omega}$.
    Then it follows that  
    $u\in\mathcal{D}^{0,q}(\Omegaa)$. 
    Since we can write
    \begin{eqnarray*}
      u=\zeta u=\zeta\dbar N_{q-1}\dbarstar u + \zeta\dbarstar N_{q+1}\dbar 
      u,	
    \end{eqnarray*}
    we obtain the estimate
    \begin{eqnarray*}
      \|u\|^{2}
      \lesssim
      \|\zeta\dbar N_{q-1}\dbarstar u\|^{2}
      +\|\dbarstar N_{q+1}\dbar u\|^{2}.
    \end{eqnarray*}
    Because $\dbar u$ is a $\dbar$-closed $(0,q+1)$-form supported in 
    $W\ssubset V$, we can use  inequality (\ref{E:dbarstarNq+1})
    to estimate the last term in the above inequality, i.e.
    \begin{eqnarray}\label{E:firstEonu}
      \|u\|^{2}
      \lesssim
      \|\zeta\dbar N_{q-1}\dbarstar u\|^{2} +
      \delta^{2\epsilon}\|\dbar u\|^{2}+\delta^{-2+2\epsilon}\|\dbar 
      u\|_{-1}^{2},
    \end{eqnarray}
    So we are left with estimating $\|\zeta\dbar N_{q-1}\dbarstar 
    u\|^{2}$:
    \begin{eqnarray*}
     \|\zeta\dbar N_{q-1}\dbarstar u\|^{2}
      &=&
      (\lbrack\zeta^{2},\dbar\rbrack N_{q-1}\dbarstar u,
      \dbar N_{q-1}\dbarstar u)
      +(\dbar\zeta^{2}N_{q-1}\dbarstar u,
      \dbar N_{q-1}\dbarstar u)\\
      &=&
      (\lbrack\zeta^{2},\dbar\rbrack N_{q-1}\dbarstar u,u-
      \dbarstar N_{q+1}\dbar u)
      +(\dbar\zeta^{2} N_{q-1}\dbarstar u,N_{q}\dbar\dbarstar u),
    \end{eqnarray*}
    since $\dbar N_{q-1}\dbarstar u=N_{q}\dbar\dbarstar u$ for $u\in
    \mathcal{D}^{0,q}(\Omegaa)$. By our choice of the cut-off function
    $\zeta$ it follows, that the supports of 
    $\lbrack\zeta^{2},\dbar\rbrack N_{q-1}\dbarstar u$ and $u$ are disjoint. Therefore
     \begin{eqnarray*}
      \|\zeta\dbar N_{q-1}\dbarstar u\|^{2}
      \lesssim
      \underbrace{\|\lbrack\zeta^{2},\dbar\rbrack \dbarstar N_{q} u\|
      \,\|\dbarstar N_{q+1}\dbar u\|}_{(A)}
      +\underbrace{\|\dbarstar N_{q}(\dbar\zeta^{2}N_{q-1}\dbarstar u)\|\,
      \|\dbarstar u\|}_{(B)}.
    \end{eqnarray*}
    Using the (sc)-(lc) inequality, we get
    \begin{eqnarray*}
    (A)\lesssim \eta\|\lbrack\zeta^{2},\dbar\rbrack\dbarstar N_{q} u\|^{2}
    +\frac{1}{\eta}\|\dbarstar N_{q+1}\dbar u\|^{2}
    \end{eqnarray*}
    for $\eta>0$.
    Recall that $\dbarstar N_{q+1}$ is a bounded map from $L_{(0,q+1)}^{2}
    (\Omegaa)$ to $L_{(0,q)}^{2}(\Omegaa)$, and also note that 
    $\lbrack\zeta^{2},\dbar\rbrack$ is a differential operator of 
    order zero. Using inequality (\ref{E:dbarstarNq+1}) again,
    we obtain
    \begin{eqnarray*}
      (A)\lesssim
      \eta\|u\|^{2}+\frac{1}{\eta}(\delta^{2\epsilon}\|\dbar u\|^{2}
      +\delta^{-2+2\epsilon}\|\dbar u\|_{-1}^{2}).
    \end{eqnarray*}	
    To estimate term $(B)$ note that $\dbar\zeta^{2}N_{q-1}\dbarstar u$ is a 
    $\dbar$-closed $(0,q)$-form, which is supported in $V$. Thus, by 
    our estimate (\ref{E:dbarstarNq}) on $\dbarstar N_{q}$ , 
    it follows 
    \begin{eqnarray*}	
      \|\dbarstar N_{q}(\dbar\zeta^{2} N_{q-1}\dbarstar 
      u)\|
      \lesssim
      \delta^{\epsilon}\underbrace{\|\dbar\zeta^{2} 
      N_{q-1}\dbarstar u\|}_{(B_{1})}
      +\delta^{-1+\epsilon}\underbrace{\|\dbar\zeta^{2} N_{q-1}\dbarstar 
      u\|_{-1}}_{(B_{2})}.
    \end{eqnarray*}
    By commuting $\dbar$ and $\zeta^{2}$, we obtain for $(B_{1})$:
    \begin{eqnarray*}
      (B_{1})
      \leq
      \|\zeta^{2}\dbar\dbarstar N_{q} u\|
      +\|\lbrack\zeta^{2},\dbar\rbrack \dbarstar N_{q} u\|
      \lesssim
      \|\dbar\dbarstar N_{q} u\|
      +\|\dbarstar N_{q} u\|
      \lesssim
      \|u\|,
    \end{eqnarray*}
    The last step holds, since $\dbar\dbarstar N_{q}$ is 
    a bounded operator on $L_{0,q}^{2}(\Omegaa)$ and 
    $\dbarstar N_{q}$ is a bounded operator from $L_{0,q}^{2}(\Omegaa)$
    to $L_{0,q-1}^{2}(\Omegaa)$.
    
    For estimating $(B_{2})$ commute $\dbar$ and $\zeta^{2}$ again, that is
    \begin{eqnarray*}
     \begin{split}
      (B_{2})
      \leq
      \|\zeta^{2}\dbar\dbarstar N_{q}u\|_{-1}
      &+\|\lbrack\zeta^{2},\dbar\rbrack\dbarstar N_{q} u\|_{-1}\\
      &\lesssim
      \|\dbar\dbarstar N_{q}u\|_{-1}
      +\|\dbarstar N_{q} u\|_{-1}
      \lesssim
      \|u\|_{-1}.
      \end{split}
    \end{eqnarray*}
    by (\ref{E:dbarstarN-1}) and (\ref{E:dbardbarstarN-1}). Combining 
    our estimates for $(B_{1})$ and $(B_{2})$, we get 
    \begin{eqnarray*}
      (B)
      \lesssim
      (\delta^{\epsilon}\|u\|
      +\delta^{-1+\epsilon}\|u\|_{-1}\|)\|\dbarstar u\|
      \lesssim
      \eta(\|u\|^{2}
      +\delta^{-2}\|u\|_{-1}^{2})
      +\frac{\delta^{2\epsilon}}{\eta}\|\dbarstar u\|^{2},
    \end{eqnarray*}
    where the last step, again, follows by the (sc)-(lc) inequality, 
    and $\eta>0$. 
    
    Recall that we need the above estimates on $(A)$ and $(B)$ to 
    get control on the term $\|\zeta\dbar N_{q-1}\dbarstar u\|$. We 
    now have
    \begin{eqnarray*}
      \|\zeta\dbar N_{q-1}\dbarstar u\|^{2}
      \lesssim
      \frac{\delta^{2\epsilon}}{\eta}(\|\dbar u\|^{2}
      +\|\dbarstar u\|^{2}
      +\delta^{-2}\|\dbar u\|_{-1}^{2})
      +\eta(\|u\|^{2}
      +\delta^{-2}\|u\|_{-1}^{2}).
    \end{eqnarray*}
    Combining this last estimate with inequality (\ref{E:firstEonu}), 
    it follows that
    \begin{eqnarray*}
      \|u\|^{2}\lesssim	
      \frac{\delta^{2\epsilon}}{\eta}(\|\dbar u\|^{2}
      +\|\dbarstar u\|^{2}
      +\delta^{-2}\|\dbar u\|_{-1}^{2})
      +\eta(\|u\|^{2}
      +\delta^{-2}\|u\|_{-1}^{2})
    \end{eqnarray*}
    holds uniformly for all $\eta>0$. 
    Finally, for all sufficiently small $\eta>0$ we can absorb the term $\eta\|u\|^2$ 
    into the left hand side and obtain
    \begin{eqnarray*}
      \|u\|^{2}\lesssim	
      \frac{\delta^{2\epsilon}}{\eta}(\|\dbar u\|^{2}
      +\|\dbarstar u\|^{2}
      +\delta^{-2}\|\dbar u\|_{-1}^{2})
      +\eta\delta^{-2}\|u\|_{-1}^{2}.	
    \end{eqnarray*}	
    Recall that here $\|.\|$ denotes the $L^{2}$-norm on $\Omegaa$. 
    However, $\Omegaa\subset\Omega$ and 
    $u\in\mathcal{D}^{0,q}(\Omega)$ is supported in $W\cap\Omegaa$. 
    Thus we can conclude
    \begin{eqnarray*}
      \|u\|_{\Omega}^{2}
      \lesssim
      \frac{\delta^{2\epsilon}}{\eta}(\|\dbar u\|_{\Omega}^{2}
      +\|\dbarstar u\|_{\Omega}^{2}
      +\delta^{-2}\|\dbar u\|_{-1,\Omega}^{2})
      +\eta\delta^{-2}\|u\|_{-1,\Omega}^{2}.	
    \end{eqnarray*}
    for all $\eta>0$ sufficiently small.
  \end{proof}    
 
 \section{Subelliptic estimate}\label{S:sectionSubest}

  In this section we show how to derive subelliptic estimates from the 
  family of estimates obtained in Proposition \ref{P:NpassingQ}. We 
  begin with stating the main result of this section.
  \begin{theorem}\label{T:SubEstonD}
    Let $\Omega\ssubset\mathbb{C}^{n}$ be a smoothly bounded domain, 
    $p$ a point on the boundary of $\Omega$. Let 
    $V$ be a special boundary chart near $p$ such that $V\cap 
    b\Omega$ is pseudoconvex. Suppose that
    \begin{eqnarray}\label{E:deltaEstonD}
      \|u\|^{2}
      \lesssim
      \frac{\delta^{2\epsilon}}{\eta}(Q(u,u)
      +\delta^{-2}\|\dbar u\|_{-1}^{2})
      +\eta\delta^{-2}\|u\|_{-1}^{2}
    \end{eqnarray}
    holds for all $u\in\mathcal{D}^{0,q}(\Omega)$ supported
    in $V\cap\bar{\Omega}$, and for all $\eta,\,\delta>0$ sufficiently small.
    Let $W\ssubset V$ be a neighborhood of $p$. Then
    \begin{eqnarray*}
      \Tnorm u\Tnorm_{\epsilon}^{2}\lesssim
      Q(u,u)
    \end{eqnarray*}	
    holds for all $u\in\mathcal{D}^{0,q}(\Omega)$ which are supported 
    in $W\cap\bar{\Omega}$.
  \end{theorem}
  For the proof of Theorem \ref{T:SubEstonD} we use a method 
  from \cite{Cat87}. That is, we introduce a 
  sequence of pseudo-differential operators, which represent a 
  partition of unity in the tangential Fourier transform variables:
  
  Let $\lbrace p_{k}(t)\rbrace_{k=0}^{\infty}$ be a sequence of 
  functions on $\mathbb{R}$ satisfying the following conditions:
  \begin{enumerate}
    \item[(1)]  $\sum_{k=0}^{\infty}p_{k}^{2}(t)=1$ for all $t\in
      \mathbb{R}$,
    \item[(2)] $p_{0}(t)=0$ for all $t\geq 2$, and $p_{k}(t)=0$
      for all $t\notin (2^{k-1},2^{k+1})$, $k\geq 1$.
  \end{enumerate}    
  We can choose the $p_{k}$'s such that $|p_{k}'(t)|\leq
  C2^{-k}$ holds for all $k\in\mathbb{N}_{0}$, $t\in\mathbb{R}$ for 
  some $C>0$.
  Let $\mathcal{S}(\mathbb{R}^{2n})$ be the class of Schwartz 
  functions on $\mathbb{R}^{2n}$. Denote by $\mathbb{R}_{-}^{2n}$ the set 
  $\lbrace (x_{1},\dots,x_{2n-1},r)\,|\,r\leq 0\rbrace$ and
  $\mathcal{S}(\mathbb{R}_{-}^{2n})$ be the restriction of 
  $\mathcal{S}(\mathbb{R}^{2n})$ to $\mathbb{R}_{-}^{2n}$.\\
  For $f\in\mathcal{S}(\mathbb{R}_{-}^{2n})$ define the 
  operators $P_{k}$ by
  \begin{eqnarray*}
    \widetilde{P_{k}f}(\xi,r):=p_{k}(|\xi|)\tilde{f}(\xi,r),  
  \end{eqnarray*}
  where $\tilde{f}$ is the tangential Fourier transform, that is
  \begin{eqnarray*}
    \tilde{f}(\xi,r)=\int_{\mathbb{R}^{2n-1}} e^{-2\pi i\langle x,\xi\rangle}
    f(x,r)dx.
  \end{eqnarray*}
  On $(0,q)$-forms we define the $P_{k}$'s to act componentwise.
  
  One of the crucial features of such operators $P_{k}$ is 
  that it makes the tangential Sobolev $s$-norm of a function $f\in 
  \mathcal{S}(\mathbb{R}_{-}^{2n})$ comparable to a series 
  involving $L^{2}$-norms 
  of $P_{k}f$. In general, we have:
  \begin{lemma}\label{L:PkTangsComp}
    For $f\in \mathcal{S}(\mathbb{R}_{-}^{2n})$ and 
    $s=s_{1}+s_{2}$ it holds that
    \begin{eqnarray*}
      \Tnorm f\Tnorm_{s}^{2}\cong\sum_{k=0}^{\infty} 
      2^{2ks_{1}}\Tnorm P_{k}f\Tnorm_{s_{2}}^{2}.	
    \end{eqnarray*}
  \end{lemma}    
  \begin{proof}
    Let $f\in \mathcal{S}(\mathbb{R}_{-}^{2n})$, 
    $s=s_{1}+s_{2}$. From the definition of the 
    tangential Sobolev $s$-norm and since $\sum_{k=0}^{\infty}p_{k}^{2}=1$ holds, it follows that 
    \begin{eqnarray*}
      \Tnorm f\Tnorm_{s}^{2}
     =
      \int_{-\infty}^{0}\int_{\mathbb{R}^{2n-1}}
      (1+|\xi|^{2})^{s}
      (\sum_{k=0}^{\infty}p_{k}^{2}(|\xi|))
      |\tilde{f}(\xi,r)|^{2}d\xi dr.
    \end{eqnarray*}  
  Since $(1+|\xi|^{2})^{s_{1}}\cong 2^{2ks_{1}}$ as long as 
    $|\xi|$ is in the support of $p_{k}$, we obtain
    \begin{eqnarray*}
      \Tnorm f\Tnorm_{s}^{2} 	
      &\cong &
      \sum_{k=0}^{\infty} 2^{2ks_{1}}\int_{-\infty}^{0}
      \int_{\mathbb{R}^{2n-1}}
      (1+|\xi|^{2})^{s_{2}}
      |p_{k}(|\xi|)\tilde{f}(\xi,r)|^{2}d\xi dr\\
      &=&
      \sum_{k=0}^{\infty} 2^{2ks_{1}}\Tnorm P_{k}f\Tnorm_{s_{2}}^{2}.
    \end{eqnarray*}	
  \end{proof}  
  Suppose $u=\psum_{|J|=q}v_Jd\bar{z}^J$ is in $\mathcal{D}^{0,q}(\Omega)$ and supported in  $V\cap\bar{\Omega}$, where $V$ is a  
  special boundary chart near a boundary point $p$. Then we can write 
  \begin{eqnarray*}
    u=\psum_{|I|=q}u_I dx^{I},
 \end{eqnarray*}  
 where $I=\{i_1,...,i_q\}$ with $1\leq i_l\leq 2n$. The operator  $P_k$ acting on a $(0,q)$-form $u$ means the following:
 \begin{eqnarray*}
 P_k u=\psum_{|I|=q}(P_k u_I)dx^I.
 \end{eqnarray*}
 We remark that $u\in\mathcal{D}^{0,q}(\Omega)$ if and only if 
 $u_I(x',0)=0$ for $x'\in\mathbb{R}^{2n-1}$ whenever $2n\in I$.
  This leads to another crucial property of the operator $P_{k}$, that is:
  $P_{k}u\in\mathcal{D}^{0,q}(\Omega)$ whenever $u\in
  \mathcal{D}^{0,q}(\Omega)$. However, the $P_{k}$'s do not see 
  the support of $u$, i.e. if $u$ is compactly supported, we can not 
  conclude the same for $P_{k}u$. Thus inequality 
  (\ref{E:deltaEstonD}) does not hold for $P_{k}u$ in general. 
  We shall introduce an appropriately chosen cut-off 
  function $\chi$ and consider $\chi P_{k}u$. To be able to deal with 
  certain error terms arising from inequality (\ref{E:deltaEstonD})
  applied to $\chi P_{k}u$, 
  we collect a few facts in the following lemmata.
  \begin{lemma}\label{L:PkfuncComE}
    If $f,g\in\mathcal{S}(\mathbb{R}_{-}^{2n})$ and
    $\sigma\in\mathbb{R}$, then 
    \begin{eqnarray*}
      \sum_{k=0}^{\infty} 2^{2k\sigma}\|\lbrack P_{k},f\rbrack g\|^{2}
      \lesssim\Tnorm g\Tnorm_{\sigma-1}^{2},
    \end{eqnarray*}	
    where the constant in $\lesssim$ does not depend on $g$.
  \end{lemma}    
   The proof of Lemma \ref{L:PkfuncComE} 
    \begin{lemma}\label{L:PkDiff-1E}
    Let $D$ be any differential 
    operator of first order 
    with coefficients in $C^{\infty}(\mathbb{R}_{-}^{2n})$ acting on 
    smooth $q$-forms, let 
    $\chi\in\mathcal{S}(\mathbb{R}_{-}^{2n})$ and $\sigma>0$. Then 
    \begin{eqnarray*}
      \sum_{k=0}^{\infty}2^{2k\sigma}\|D(\chi P_{k}u)\|_{-\sigma}^{2}
      \lesssim \|Du\|^{2}+\|u\|^{2}+\psum_{|I|=q}\Tnorm\frac{\partial 
      u_{I}}{\partial x_{2n}}\Tnorm_{-1}^{2}.
    \end{eqnarray*}
    holds for all $q$-forms $u$ with coefficients in 
    $\mathcal{S}(\mathbb{R}_{-}^{2n})$. Here, the constant in 
    $\lesssim$ does not depend on $u$.
  \end{lemma}  
    \begin{proof}
       Recall that $\Lambda_{t}^{-\sigma}$ denotes the tangential 
       Bessel potential of order $-\sigma$. We obtain
       \begin{eqnarray*}
	 \sum_{k=0}^{\infty}2^{2k\sigma}\|D(\chi P_{k}u)\|_{-\sigma}^{2}
	 &\leq&
	 \sum_{k=0}^{\infty}2^{2k\sigma}\Tnorm 
	 D(\chi P_{k}u)\Tnorm_{-\sigma}^{2}
	 =\sum_{k=0}^{\infty}2^{2k\sigma}\|\Lambda_{t}^{-\sigma}
	 D(\chi P_{k}u)\|^{2}\\
	 &\lesssim&
	 \sum_{k=0}^{\infty}2^{2k\sigma}\|\chi\Lambda_{t}^{-\sigma}
	 D P_{k}u\|^{2}+
	 \sum_{k=0}^{\infty}2^{2k\sigma}\|\lbrack\Lambda_{t}^{-\sigma}
	 D,\chi\rbrack P_{k}u\|^{2},
       \end{eqnarray*}	    
       where the last step follows by commuting.
       We note that $\lbrack\Lambda_{t}^{-\sigma}D,\chi\rbrack$ is of 
       tangential order $-\sigma$ and of normal order $0$. Therefore, 
       invoking Lemma \ref{L:PkTangsComp}, we get
       \begin{eqnarray*}
	 \sum_{k=0}^{\infty}2^{2k\sigma}\|\lbrack\Lambda_{t}^{-\sigma}
	 D,\chi\rbrack P_{k}u\|^{2}
	 \lesssim\sum_{k=0}^{\infty}2^{2k\sigma}\Tnorm P_{k}u
	 \Tnorm_{-\sigma}^{2}
	 \cong
	 \|u\|^{2}.
       \end{eqnarray*}
       Similarly, we obtain by commuting
       \begin{equation*}
	   \begin{split}
	 \sum_{k=0}^{\infty}2^{2k\sigma}&\|\chi\Lambda_{t}^{-\sigma}
	 D P_{k}u\|^{2}
	 \lesssim
	  \sum_{k=0}^{\infty}2^{2k\sigma}\|\chi D\Lambda_{t}^{-\sigma}
	 P_{k}u\|^{2}+\sum_{k=0}^{\infty}2^{2k\sigma}\Tnorm 
	 P_{k}u\Tnorm_{-\sigma}^{2}\\
	 &\lesssim
	 \sum_{k=0}^{\infty}2^{2k\sigma}\|P_{k}
	 (\chi D\Lambda_{t}^{-\sigma}
	 u)\|^{2}+
	 \sum_{k=0}^{\infty}2^{2k\sigma}\|\lbrack 
	 \chi D\Lambda_{t}^{-\sigma},
	 P_{k}\rbrack u\|^{2}+\|u\|^{2}\\
	 &\lesssim
	 \underbrace{\Tnorm
	 \chi D\Lambda_{t}^{-\sigma}
	 u\Tnorm_{\sigma}^{2}}_{(A)}+
	 \sum_{k=0}^{\infty}2^{2k\sigma}\underbrace{\|\lbrack 
	 \chi D\Lambda_{t}^{-\sigma},
	 P_{k}\rbrack u\|^{2}}_{(B_{k})}+\|u\|^{2},
	 \end{split}
       \end{equation*}
       where the last line follows again by Lemma \ref{L:PkTangsComp}.
       We write
       \begin{eqnarray*}
	 \chi D=\psum_{|I|=q}\sum_{j=1}^{2n}a_{j}^{I}\frac{\partial}{\partial 
	 x_{j}},  
       \end{eqnarray*}	   
       and estimate term (A) by commuting:
       \begin{eqnarray*}
	 (A)&=&\Tnorm\chi D\Lambda_{t}^{-\sigma} u\Tnorm_{\sigma}^{2}
	 \lesssim
	 \Tnorm\Lambda_{t}^{-\sigma}\chi D u\Tnorm_{\sigma}^{2}
	 +\Tnorm\lbrack\chi D,\Lambda_{t}^{-\sigma}\rbrack 
	 u\Tnorm_{\sigma}^{2}\\
	 &\lesssim&
	 \|Du\|^{2}+\psum_{|I|=q}\sum_{j=1}^{2n}
	 \Tnorm \lbrack a_{j}^{I}\frac{\partial}{\partial x_{j}}
	 ,\Lambda_{t}^{-\sigma}\rbrack u_{I}\Tnorm_{\sigma}^{2}.
       \end{eqnarray*}
       Since $\frac{\partial}{\partial x_{j}}$ and 
       $\Lambda_{t}^{-\sigma}$ commute, it follows that
       \begin{eqnarray*}
         (A)&\lesssim&\|Du\|^2+\psum_{|I|=q}\sum_{j=1}^{2n}
	 \Tnorm \lbrack a_{j}^{I}
	 ,\Lambda_{t}^{-\sigma}\rbrack \frac{\partial u_{I}}
	 {\partial x_{j}}\Tnorm_{\sigma}^{2}\\
	 &\lesssim&
	 \|Du\|^2+
	 \|u\|^{2}+\psum_{|I|=q}\Tnorm\frac{\partial u_{I}}
	 {\partial x_{2n}}\Tnorm_{-1}^{2}.
       \end{eqnarray*}
       Here, the last estimate holds since $\lbrack 
       a_{j}^{I},\Lambda_{t}^{-\sigma}\rbrack$ is of tangential 
       order $-\sigma-1$ and $\frac{\partial}{\partial x_{j}}$ is a 
       tangential derivative if $j\in\{1,\ldots,2n-1\}$. We are left 
       with estimating the terms $(B_{k})$. We first notice that
       \begin{eqnarray*}
	 (B_{k})
	 \lesssim
	 \psum_{|I|=q}\sum_{j=1}^{2n}\|\lbrack a_{j}^{I},P_{k}\rbrack 
	 \frac{\partial}{\partial x_{j}}\Lambda_{t}^{-\sigma} u_{I}\|^{2}.
       \end{eqnarray*}	   
       Lemma \ref{L:PkfuncComE} implies now
       \begin{eqnarray*}
	 \sum_{k=0}^{\infty}2^{2k\sigma}(B_{k})
	 \lesssim
	 \psum_{|I|=q}\sum_{j=1}^{2n}\Tnorm\frac{\partial}{\partial x_{j}}
	 \Lambda_{t}^{-\sigma}u_{I}\Tnorm_{\sigma-1}^{2}
	 \lesssim
	 \|u\|^{2}+\psum_{|I|=q}\Tnorm\frac{\partial u_{I}}{\partial x_{2n}}
	 \Tnorm_{-1}^{2}.
       \end{eqnarray*}	 
       Combining all our estimates we end up with the claimed 
       inequality.
       \begin{eqnarray*}
	 \sum_{k=0}^{\infty}2^{2k\sigma}\|D(\chi P_{k}u)\|_{-\sigma}^{2}
         \lesssim \|Du\|^{2}+\|u\|^{2}+\psum_{|I|=q}\Tnorm\frac{\partial 
         u_{I}}{\partial x_{2n}}\Tnorm_{-1}^{2}.  
       \end{eqnarray*}	   
    \end{proof}
    Having collected the basic facts concerning the $P_{k}$'s, we are 
    ready to prove Theorem \ref{T:SubEstonD}.
    \begin{proof}[Proof of Theorem \ref{T:SubEstonD}]
	Let $V$ be a special boundary chart near p such that inequality
	(\ref{E:deltaEstonD}) holds, that is
	\begin{eqnarray*}
          \|u\|^{2}
          \lesssim
          \frac{\delta^{2\epsilon}}{\eta}(Q(u,u)
	  +\delta^{-2}\|\dbar u\|_{-1}^{2})
          +\eta\delta^{-2}\|u\|_{-1}^{2}
        \end{eqnarray*}
	holds for all $u\in\mathcal{D}^{0,q}(\Omega)$ supported
        in $V\cap\bar{\Omega}$. Let $W\ssubset V$ be a neighborhood 
        of $p$, and $u\in\mathcal{D}^{0,q}(\Omega)$ supported in 
        $W\cap\bar{\Omega}$. Let $\chi\in 
        C_{c}^{\infty}(V)$ such that $\chi=1$ on $W$ and $\chi\geq 0$. Then it follows 
        by Lemma \ref{L:PkTangsComp} and by commuting
	\begin{eqnarray*}
	  \Tnorm u\Tnorm_{\epsilon}^{2}
	  =
	  \Tnorm\chi u\Tnorm_{\epsilon}^{2}
	  &\lesssim&
	  \sum_{k=0}^{\infty}2^{2k\epsilon}\|\chi P_{k}u\|^{2}
	  +\sum_{k=0}^{\infty}2^{2k\epsilon}\|\lbrack P_{k},\chi\rbrack 
	  u\|^{2}\\
	  &\lesssim&
	  \sum_{k=0}^{\infty}2^{2k\epsilon}\|\chi P_{k}u\|^{2}
	  +\Tnorm u\Tnorm_{\epsilon-1}^{2},
	\end{eqnarray*}
	where the last step follows by Lemma \ref{L:PkfuncComE}.
	Since $\epsilon\leq\frac{1}{2}$ holds, we obtain
	\begin{eqnarray*}
	  \Tnorm u\Tnorm_{\epsilon}^{2}
	  \lesssim
	  \sum_{k=0}^{\infty}2^{2k\epsilon}\|\chi P_{k}u\|^{2}
	  +\|u\|^{2}.
	\end{eqnarray*}
	Now inequality (\ref{E:deltaEstonD}) comes into play. Since
	$\chi P_{k}u\in\mathcal{D}^{0,q}(\Omega)$ is supported in 
	$V\cap\bar{\Omega}$, it follows that
	\begin{eqnarray*}
	  \|\chi P_{k}u\|^{2}\lesssim\frac{\delta^{2\epsilon}}{\eta}
	  (Q(\chi P_{k}u,\chi P_{k}u)+
	  \delta^{-2}\|\dbar(\chi P_{k} u)\|_{-1}^{2})
	  +\eta\delta^{-2}\|\chi P_{k} u\|_{-1}^{2}
	\end{eqnarray*}
	holds uniformly for all $k\in\mathbb{N}_{0}$, for all positive 
	$\delta<\delta_{0}$ and $\eta<\eta_{0}$. 
	Let $k_{0}\in\mathbb{N}$ such that $2^{-k_{0}}\leq\delta_{0}$. Then we 
	obtain for all $k\geq k_{0}$
	\begin{eqnarray*}
	  2^{2k\epsilon}\|\chi P_{k}u\|^{2}
	  \lesssim
	  \frac{1}{\eta}(Q(\chi P_{k} u,\chi P_{k}u)+
	  2^{2k}\|\dbar(\chi P_{k} u)\|_{-1}^{2})
	  +\eta 2^{2k(1+\epsilon)}\|\chi P_{k} u\|_{-1}^{2}.
	\end{eqnarray*}
	Observe that
	\begin{eqnarray*}
	  \sum_{k=0}^{k_{0}-1}2^{2k\epsilon}\|\chi P_{k}u\|^{2}
	  \leq
	  \sum_{k=0}^{k_{0}-1}2^{2k\epsilon}\|u\|^{2}
	  \lesssim\|u\|^{2}.
	\end{eqnarray*}    
	Thus we can sum up over $k\in\mathbb{N}_{0}$, obtaining
	\begin{eqnarray*}
	  \sum_{k=0}^{\infty}2^{2k\epsilon}\|\chi P_{k}u\|^{2}
	  &\lesssim&
	  \frac{1}{\eta}
	  \sum_{k=0}^{\infty}(Q(\chi P_{k}u,\chi P_{k}u))
	  +\frac{1}{\eta}
	  \sum_{k=0}^{\infty}2^{2k}
	  \|\dbar(\chi P_{k}u)\|_{-1}^{2}\\
	  & &+\eta\sum_{k=0}^{\infty}2^{2k(1+\epsilon)}
	  \|\chi P_{k}u\|_{-1}^{2}+\|u\|^{2}
	\end{eqnarray*}
	Using Lemma \ref{L:PkTangsComp}, we have
	\begin{eqnarray*}
	  \sum_{k=0}^{\infty} 2^{2k(1+\epsilon)}\|\chi P_{k}u\|_{-1}^{2}
	  \lesssim
	  \sum_{k=0}^{\infty} 2^{2k(1+\epsilon)}\Tnorm P_{k}u\Tnorm_{-1}^{2}
	  \cong\Tnorm u\Tnorm_{\epsilon}^{2}.
	\end{eqnarray*}   
	Furthermore, applying Lemma \ref{L:PkDiff-1E} with $\sigma=0$ and 
	$\sigma=1$ resp., we get
	\begin{eqnarray*}
	  \sum_{k=0}^{\infty}Q(\chi P_{k}u,\chi P_{k}u)
	  +
	  \sum_{k=0}^{\infty}2^{2k}
	  \|\dbar(\chi P_{k}u)\|_{-1}^{2}
	  \lesssim
	  Q(u,u)+\|u\|^{2}+\psum_{|I|=q}\Tnorm \frac{\partial u_{I}}{\partial 
	  x_{2n}}\Tnorm_{-1}^{2}.
	\end{eqnarray*}    
	Note that $\frac{\partial}{\partial x_{2n}}$ can be expressed as a 
	linear combination of the $\frac{\partial}{\partial\bar{z}_{j}}$'s 
	and a tangential vector field $T$.
	Then
	\begin{eqnarray*}
	  \Tnorm\frac{\partial u_{I}}{\partial x_{2n}}\Tnorm_{-1}^{2}
	  &\lesssim&
	  \sum_{j=1}^{n}\Tnorm\frac{\partial u_{I}}{\partial\bar{z}_{j}}
	  \Tnorm_{-1}^{2}
	  +\Tnorm Tu_{I}\Tnorm_{-1}^{2}
	  \lesssim
	  \sum_{j=1}^{n}\|\frac{\partial u_{I}}{\partial\bar{z}_{j}}\|^{2}
	  +\|u_{I}\|^{2}\\
	  &\lesssim&
	  \|\dbar u\|^{2}+\|\dbarstar u\|^{2}+\|u\|^{2}
	  \lesssim Q(u,u).
	\end{eqnarray*}    
	Thus, by combining our estimates , we obtain
	\begin{eqnarray*}
	  \Tnorm u\Tnorm_{\epsilon}^{2}
	  \lesssim
	  \sum_{k=0}^{\infty}2^{2k\epsilon}\|\chi P_{k}u\|^{2}
	  +\|u\|^{2}
	  \lesssim
	  \frac{1}{\eta}Q(u,u)+\eta\Tnorm u\Tnorm_{\epsilon}^{2}.
	\end{eqnarray*}
	Choosing $\eta>0$ small enough, we can absorb the term
	$\eta\Tnorm u\Tnorm_{\epsilon}^{2}$ into the left hand side and it follows
	 $ \Tnorm u\Tnorm_{\epsilon}^{2}\lesssim Q(u,u)$. 
    \end{proof}	    

\section{An Example}
Consider the domain 
$D=\{w\in\mathbb{C}^3\;|\;\rho(w):=\text{Re}\;w_3+|w_1^2-w_2 w_3|^2+|w_2^2|^2<0\}$ near the origin. The  $1$-type (in the sense of D'Angelo \cite{D'An82}) at  $(0,0,0)$ is $4$, but at any boundary point of the form $(0,0,i\epsilon)$, $\epsilon>0$, the $1$-type is $8$. In the following we show that a subelliptic estimate of order $\frac{1}{8}-\eta$ holds for any $\eta>0$ near the origin.  Instead of constructing the $\{\phidel\}$ on $D$, we 
consider $\Omega=\{z\in\mathbb{C}^3\;|\;r(z)<0\}$, where
\begin{eqnarray*}
  r(z)=|z_3|^2-1+|(1+z_3)z_1^2-z_2(z_3-1)|^2+|1+z_3|^2|z_2|^4<0,
\end{eqnarray*}
in a neighborhood $U$ of the boundary point $p=(0,0,1)$. 
Notice that $D$ near the origin is biholomorphic to $\Omega$ via the transformation $z_1=w_1$, $z_2=w_2$ and $z_3=\frac{w_{3}+1}{1-w_3}$. We claim that 
\begin{eqnarray*}
\phidel(z)&=&-\log(-r(z)+\delta)-\log(-\log(|z_1|^2+\delta^\frac{1}{4}))\\
& &-\log(-\log(|z_2|^2+\delta^{\frac{1}{2}+\eta}))\\
&=&\psi_0(z)+\psi_1(z)+\psi_2(z)
\end{eqnarray*}
satisfies the hypotheses of Theorem \ref{T:MainTheorem} on $\overline{\Omega}\cap U$ with
$\epsilon=\frac{1}{8}-\frac{\eta}{2}$ for $\eta>0$. A straightforward computation shows that $\phidel$ is plurisubharmonic and has a self-bounded complex gradient near $p$. In the following we show that
\begin{eqnarray}\label{E:ExampleHess}
i\partial\dbar\phidel(z)(\xi,\xi)\geq C\delta^{-\frac{1}{4}+\eta}|\xi|^2
\end{eqnarray}
holds for all $\xi\in\mathbb{C}^3$ and $z\in\esdel\cap U$. One computes 
\begin{eqnarray*}
  i\partial\dbar r(\xi,\xi)&=&4|1+z_3|^2|z_1|^2|\xi_1|^2 +(|z_3-1|^2+4|z_2|^2|1+z_3|^2)|\xi_2|^2\\
  & &+(1+|z_1^2-z_2|^2+|z_2|^4)|\xi_3|^2\\
  & &+2\text{Re}((2(1+z_3)^2z_2 \bar{z}_2^2-(z_3-1)(\bar{z}_1^2-\bar{z}_2))\xi_2\bar{\xi}_3)\\
  & &+4\text{Re}((1+z_3)z_1\xi_1((\bar{z}_1^2-\bar{z}_2)\bar{\xi}_3-
  (\bar{z}_3-1)\bar{\xi}_2).
\end{eqnarray*}
Denote the last term on the right hand side by (I). Estimating (I) we obtain
\begin{eqnarray*}
(I)&\geq& -4|1+z_3|^2|z_1|^2|\xi_1|^2-|z_3-1|^2|\xi_2|^2-|z_1^2-z_2|^2|\xi_3|^2\\
& &+2\text{Re}((z_3-1)\xi_2(\bar{z}_1^2-\bar{z}_2)\bar{\xi}_3).
\end{eqnarray*}
It follows easily that
\begin{eqnarray}\label{E:Hessr1}
i\partial\dbar r(z)(\xi,\xi)\geq|z_2|^2|\xi_2|^2+\frac{1}{2}|\xi_3|^2.
\end{eqnarray}
This estimate implies that if $z\in\esdel\cap U$, then
\begin{eqnarray*}
  i\partial\dbar \psi_0(z)(\xi,\xi)
  \geq\frac{|z_2|^2|\xi_2|^2+\frac{1}{2}|\xi_3|^2}{-r(z)+\delta}
  \geq\frac{1}{4}(\delta^{-\frac{1}{2}}|\xi_2|^2+\delta^{-1}|\xi_3|^2),
\end{eqnarray*}
where the first estimate on the right hand side only holds if $|z_2|^2\geq\delta^{\frac{1}{2}}$. If $|z_2|^2\leq\delta^{\frac{1}{2}}$, then
\begin{eqnarray*}
  i\partial\dbar\psi_2(z)(\xi,\xi)\geq\frac{\delta^{\frac{1}{2}+\eta}|\xi_2|^2}
  {-\log(|z_2|^2+\delta^{\frac{1}{2}+\eta})(|z_2|^2+\delta^{\frac{1}{2}+\eta})^2}
  \geq\delta^{-\frac{1}{4}}|\xi_2|^2.
\end{eqnarray*}
Similarly, we obtain $i\partial\dbar\psi_1(z)(\xi,\xi)\geq\frac{1}{4}\delta^{-\frac{1}{4}+\eta}|\xi_1|^2$ for $|z_1|^2\leq\delta^{\frac{1}{4}}$ for all $\delta>0$ sufficiently small. Thus it remains to show that (\ref{E:ExampleHess}) holds also in the directions involving $\xi_{1}$ for $z\in\esdel\cap U$ with $|z_1|^2\geq\delta^{\frac{1}{4}}$. For that we shall use a different estimate for the complex Hessian of $r$, that is
\begin{eqnarray}\label{E:Hessr2}
  i\partial\dbar r(z)(\xi,\xi)\geq \frac{1}{2}|z_1|^2|\xi_1|^2-4|z_3-1|^2|\xi_2|^2.
\end{eqnarray}
Then, if $z\in\esdel$ and $|z_1|^2\geq\delta^{\frac{1}{4}}$, we obtain by using (\ref{E:Hessr2}) and (\ref{E:Hessr1})
\begin{eqnarray*}
  (\delta^{\frac{1}{2}+\eta}+\frac{1}{2})i\partial\dbar\psi_0(z)(\xi,\xi)
  \geq \frac{C}{\delta}(\delta^{\frac{3}{4}+\eta}|\xi_1|^2+(|z_2|^2-16\delta^{\frac{1}{2}+\eta}|z_3-1|^2)|\xi_2|^2).
\end{eqnarray*}
Thus we obtain (\ref{E:ExampleHess}) for all $z\in\esdel$ as long as 
$|z_2|^2\geq 16\delta^{\frac{1}{2}+\eta}|z_3-1|^2$. 
If the latter inequality is not true, then we can assume that $|z_2|^2\leq\delta^{\frac{1}{2}+\eta}$. However, in that case 
\begin{eqnarray*}
  \frac{1}{2}i\partial\dbar\psi_2(z)(\xi,\xi)-16\delta^{-\frac{1}{2}+\eta}|z_3-1|^2|\xi_2|^2\geq 0,
\end{eqnarray*}
which completes the proof of (\ref{E:ExampleHess}).

With a construction similar to the above one obtains for the domains 
\begin{eqnarray*}
 D_{k,l,m,n}=\{w\in\mathbb{C}^{3}\;|\;\text{Re}w_3+|w_1^k-w_2^l w_3^m|^2+|w_2^n|^2<0\},\;\;k,l,m,n\in\mathbb{N}
\end{eqnarray*}
a subelliptic estimate of  order $\frac{1}{M}-\eta$, $\eta>0$, where $M$ is the maximum $1$-type near the origin.

\end{document}